\documentclass[12pt]{article}
\usepackage{amsmath, amsfonts, amssymb, amsthm, color}
\usepackage[hypertex,draft=false]{hyperref}

\newtheorem{theorem}{Theorem}
\newtheorem{proposition}[theorem]{Proposition}
\newtheorem{lemma}[theorem]{Lemma}
\newtheorem{corollary}[theorem]{Corollary}
\theoremstyle{definition}
\newtheorem{definition}[theorem]{Definition}
\newtheorem{definition and remark}[theorem]{Definition and Remark}
\newtheorem{example}[theorem]{Example}
\newtheorem{remark}[theorem]{Remark}
\newtheorem{notation}[theorem]{Notation}
\newtheorem{notation and remark}[theorem]{Notation and Remark}

\DeclareMathOperator{\rank}{rank}
\DeclareMathOperator{\Ker}{Ker}
\DeclareMathOperator{\Tor}{Tor}
\DeclareMathOperator{\In}{in}
\DeclareMathOperator{\gin}{gin}

\newcommand{\frakm}{\mathfrak{m}}
\newcommand{\frakp}{\mathfrak{p}}
\newcommand{\inrevlex}{\In_{\text{revlex}}}
\newcommand{\PPii}{K[x_1,x_2]}
\newcommand{\PPiii}{K[x_1,x_2,x_3]}
\newcommand{\PPiv}{K[x_1,x_2,x_3,x_4]}
\newcommand{\PPn}{K[x_1,x_2,\ldots,x_n]}
\newcommand{\PPnm}{K[x_1,x_2,\ldots,x_{n-1}]}
\newcommand{\Rbar}{\overline{R}}
\newcommand{\Ibar}{\overline{I}}

\newcommand{\HH}{{\mathbf H}}
\newcommand{\cM}{{\mathcal M}}

\newcommand{\cB}{{\mathcal B}}
\newcommand{\cA}{{\mathcal A}}

\newcommand{\NOTYET}[1]%
{{\color{red}$\blacksquare$}\marginpar{\scriptsize\sf{
      \makebox[0pt][l]{\parbox[t]{22ex}{#1}}}}}

\title{
  Generic initial ideals, graded Betti numbers
  and $k$-Lefschetz properties
}
\author{Tadahito Harima and Akihito Wachi}
\begin{document}

\def\thefootnote{\relax}
\footnotetext{\hspace{-5ex}
}

\date{}

\maketitle

Deditcated to Professor Junzo Watanabe on his sixtieth birthday

\begin{abstract}
We introduce the $k$-strong Lefschetz property ($k$-SLP) 
and the $k$-weak Lefschetz property ($k$-WLP) 
for graded Artinian $K$-algebras,
which are generalizations of the Lefschetz properties.
The main results obtained in this paper are as follows:

1. Let $I$ be a graded ideal of $R=K[x_1, x_2, x_3]$ 
whose quotient ring $R/I$ has the SLP.
Then the generic initial ideal of $I$ is the unique
almost revlex ideal with the same Hilbert function as $R/I$.

2. Let $I$ be a graded ideal of $R=K[x_1, x_2, \ldots, x_n]$ 
whose quotient ring $R/I$ has the $n$-SLP.
Suppose that all $k$-th differences of the Hilbert function of $R/I$
are quasi-symmetric.
Then the generic initial ideal of $I$
is the unique almost revlex ideal 
with the same Hilbert function as $R/I$.

3. We give a sharp upper bound on the graded Betti numbers 
of Artinian $K$-algebras with the $k$-WLP 
and a fixed Hilbert function.
\end{abstract}

\tableofcontents
\section{Introduction}
\label{sec:introduction}
The strong and weak Lefschetz properties for graded Artinian
$K$-algebras
(Definition \ref{defn:slp}. SLP and WLP for short)
are often used in studying
generic initial ideals and graded Betti numbers
(\cite{MR951211}, \cite{MR1184062}, \cite{MR2033004},
  \cite{MR1970804}, \cite{math-0610649},  \cite{math-0610733},
  \cite{ahn-shin-jpaa}).
We generalize the Lefschetz properties,
and define the $k$-strong Lefschetz property ($k$-SLP)
and the $k$-weak Lefschetz property ($k$-WLP) for graded Artinian
$K$-algebras
(Definition~\ref{defn:k-SLP-k-WLP}).
This notion was first introduced by A.~Iarrobino 
in a private conversation with J.~Watanabe in 1995. 
The first purpose of the paper
(Theorem~\ref{thm:gin-of-symmetric-nSLP}) 
is to determine the generic initial ideals of ideals
whose quotient rings have the $n$-SLP
and have Hilbert functions satisfying some condition.
The second purpose is 
to give upper bounds of the graded Betti numbers 
of graded Artinian $K$-algebras with the $k$-WLP
(Theorem~\ref{thm-bounds}).

Let $K$ be a field of characteristic zero.
Suppose that a graded Artinian
$K$-algebra $A$ has the SLP (resp. WLP),
and $\ell \in A$ is a Lefschetz element.
If $A/(\ell)$ again has the SLP (resp. WLP),
then we say that $A$ has the 2-SLP (resp. 2-WLP).
We recursively define the $k$-Lefschetz properties
(Definition~\ref{defn:k-SLP-k-WLP}):
$A$ is said to have the $k$-SLP (resp. $k$-WLP),
if $A$ has the SLP (resp. WLP), 
and $A/(l)$ has the $(k-1)$-SLP (resp. $(k-1)$-WLP).
We have the characterization of the Hilbert functions 
of graded Artinian 
$K$-algebras with the $k$-SLP or the $k$-WLP
(Proposition~\ref{prop:hilbert-fnct-of-nSLP}).
In addition, for a graded Artinian ideal $I \subset R$,
we show that
$R/I$ has the $k$-SLP (resp. $k$-WLP)
if and only if $R/\gin(I)$ has the $k$-SLP (resp. $k$-WLP),
where $\gin(I)$ denotes the generic initial ideal
with respect to the graded reverse lexicographic order.

We explain our results on generic initial ideals.
A monomial ideal $I$ is called an almost revlex ideal,
if the following condition holds:
for each minimal generator $u$ of $I$,
every monomial $v$ with $\deg v = \deg u$ and $v >_{\text{revlex}} u$
belongs to $I$.
Almost revlex ideals play a key role in the paper.
The characterization of the Hilbert functions for almost revlex ideals
is given in Proposition~\ref{prop:hilbert-fnct-of-almost}.
We start with the uniqueness of generic initial ideals 
in the case of three variables.
\medskip\par\noindent\textbf{%
  Theorem (see Theorem~\ref{thm:3var-gin}).
}{\it
Let $I \subset R = \PPiii$ be a graded Artinian ideal
whose quotient ring has the SLP.
Then $\gin(I)$ is the unique almost revlex ideal
for the Hilbert function of $R/I$.
In particular, $\gin(I)$ is uniquely determined 
only by the Hilbert function.
}
\medskip\par\noindent
For related results of the case of three variables,
see Cimpoea\c{s} \cite{math-0610649} and 
Ahn-Cho-Park \cite{math-0610733},
where the uniqueness of $\gin(I)$ is proved 
under slightly stronger conditions than in the theorem above.
We give some examples of complete intersection of height three
whose generic initial ideals are the unique almost revlex ideals
(Example~\ref{ex:3var-ci-SLP}).

By using the $n$-SLP,
we obtain the following result for the case of $n$ variables.
In the following theorem,
`quasi-symmetric' is a notion including `symmetric'
(Definition~\ref{defn:quasi-symmetric}).
\medskip\par\noindent\textbf{%
  Theorem (see Theorem~\ref{thm:gin-of-symmetric-nSLP}).
}{\it
  Let $I \subset \PPn$ be a graded Artinian ideal 
  whose quotient ring has the $n$-SLP,
  and has the Hilbert function $h$.
  Suppose that the $k$-th difference $\Delta^k h$ is quasi-symmetric
  for every integer $k$ with $0 \le k \le n-4$.
  Then $\gin(I)$ is the unique almost revlex ideal
  for the Hilbert function $h$.
  In particular, $\gin(I)$ is uniquely determined
  only by the Hilbert function $h$.

  Here the operator $\Delta$ is defined by
  $(\Delta h)_i = \max\{h_i - h_{i-1}, 0\}$,
  and $\Delta^k h$ is the sequence obtained by applying 
  $\Delta$ $k$-times.
}
\medskip\par\noindent
The key to proving this theorem is a uniqueness of Borel-fixed ideals
whose quotient rings have the $n$-SLP
(Theorem~\ref{thm:uniqueness-of-Borel-nSLP-symmetric}).
We give some examples of complete intersection of height $n$
whose generic initial ideals are the unique almost revlex ideals
(Example~\ref{ex:n-var-ci-SLP}).

We next explain our result on the maximality of graded Betti numbers.
Let $R = \PPn$.
The following result on the maximal graded Betti numbers
is first proved for $k=1$ by
Harima-Migliore-Nagel-Watanabe \cite{MR1970804}.
\medskip\par\noindent\textbf{%
  Theorem (see Theorem~\ref{thm-bounds} 
    and Corollary~\ref{cor:maximal-Betti-number}).
}{\it
  Let $h$ be the Hilbert function of some graded Artinian
  $K$-algebra with the $k$-WLP.
  Then there is a Borel fixed ideal $I$ of $R$ 
  such that $R/I$ has the $k$-SLP, 
  the Hilbert function of $R/I$ is $h$,
  and $\beta_{i,i+j}(A) \le \beta_{i,i+j}(R/I)$
  for all graded Artinian $K$-algebra $A$
  having the $k$-WLP and $h$ as Hilbert function,
  and for any $i$ and $j$.

  In particular, when $k=n$,
  the ideal $I$ for the upper bounds 
  is the unique almost revlex ideal for the Hilbert function $h$.
}
\medskip\par

Some of the results of this paper have been obtained
independently and at the same time by
Constantinescu (see \cite{constantinescu})
and Cho-Park (see \cite{park-0707-1365}).

\section{Generic initial ideals in $\PPiii$ and the SLP}

In this section,
we first recall the Lefschetz properties (Definition~\ref{defn:slp})
and related facts.
We then introduce the notion of $x_n$-chains
(Definition and Remark~\ref{defn:xn-chain}),
which is a useful tool to study standard monomials.
The main goal of this section is Theorem~\ref{thm:3var-gin}:
for a graded Artinian ideal $I \subset \PPiii$
whose quotient ring has the SLP,
the generic initial ideal of $I$ 
with respect to the the graded reverse lexicographic order
is the unique almost revlex ideal 
for the same Hilbert function as $\PPiii/I$.

\subsection{The Lefschetz properties and $x_n$-chain decomposition}

\begin{definition}
\label{defn:slp} 
Let $A$ be a graded Artinian algebra over a field $K$,
and $A = \bigoplus_{i=0}^c A_i$ its decomposition
into graded components.
The graded algebra $A$ is said to have the
{\it strong (resp. weak) Lefschetz property},
if there exists an element $\ell \in A_1$ such that
the multiplication map
$\times \ell^s: A_i \to A_{i+s}$ ($f \mapsto \ell^s f$)
is full-rank for every $i \ge 0$ and $s>0$ (resp. $s=1$).
In this case, $\ell$ is called a {\it Lefschetz element},
and we also say that $(A,\ell)$ has 
the strong (resp. weak) Lefschetz property.
We abbreviate these properties as the {\it SLP (resp. WLP)} for short.
\end{definition}
It is clear that 
if $(A, \ell)$ has the SLP, then $(A, \ell)$ has the WLP.
It is also clear that Hilbert functions of graded algebras
with the SLP or the WLP are {\it unimodal}.
Namely there exists a non-negative integer $i$ such that
$h_0, h_1, \ldots, h_i$ is an increasing sequence
and $h_i, h_{i+1}, \ldots$ is a weakly decreasing sequence,
where $h_j = \dim_K A_j$.
Following \cite[Remark 3.3]{MR1970804},
we define positive integers $u_1, u_2, \ldots$ as follows
for later use.
\begin{align}
\label{eq:def-of-ui}
h_0 < h_1 < \cdots < h_{u_1} = 
h_{u_1+1} = \cdots = h_{u_2-1} >
h_{u_2} = \cdots = h_{u_3-1} >
h_{u_3} = \cdots
\end{align}

Suppose that the Hilbert function 
of the graded Artinian algebra $A$ is {\it symmetric},
that is, $A = \bigoplus_{i=0}^c A_i$ ($A_c \ne (0)$)
and $\dim_{K} A_i = \dim_{K} A_{c-i}$
for $i = 0, 1, \ldots, \lfloor c/2 \rfloor$.
In this case,
it is clear that $A$ has the SLP
if and only if there exists $\ell \in A_1$ and
$\times \ell^{c-2i}: A_i \to A_{c-i}$
is bijective for every $i = 0, 1, \ldots, \lfloor c/2 \rfloor$.

For a graded algebra $A$, we denote its Hilbert function by $\HH_A$.
Namely $\HH_A(t)$ denotes the linear dimension
of the graded component $A_t$ of degree $t$.
We often identify $\HH_A$ with a finite sequence
$h = (h_0, h_1, \ldots, h_c)$.
\begin{remark}[{\cite[Lemma 2.1, Lemma 2.2]{harima-watanabe-jpaa}}]
\label{rem:WLP-and-SLP}
The Lefschetz properties can be written in terms of Hilbert functions
as follows.
Let $I \subset R = \PPn$ be a graded ideal,
and $h$ the Hilbert function of $R/I$.
\begin{itemize}
  \item[(i)]
  $(R/I, \ell)$ has the WLP if and only if
  the Hilbert function of the quotient ring $R/(I+(\ell))$ is equal to
  the difference $\Delta h$ of $h$,
  where $\Delta h$ is defined by
  \begin{align*}
    (\Delta h)_i &= \max\{h_{i} - h_{i-1}, 0\}
    \qquad (i = 0, 1, 2, \ldots),
  \end{align*}
  where $h_{-1}$ is defined as zero.
  \item[(ii)]
  $(R/I, \ell)$ has the SLP if and only if the Hilbert function of
  the quotient ring $R/(I+(\ell^s))$ is equal to the sequence 
  \begin{align*}
    ( \max\{ h_i - h_{i-s} \}, 0 )_{i=0,1,\ldots},
  \end{align*}
  for every positive integer $s$.
  Here $h_i=0$ for $i<0$.
\end{itemize}
\end{remark}
A sequence
$h = (h_0, h_1, \ldots, h_c)$ is called an {\it O-sequence}
if $h$ is a Hilbert function of some graded $K$-algebra.
There is a classification of Hilbert functions
of graded Artinian algebras with the SLP or the WLP.
\begin{proposition}[{\cite[Corollary~4.6]{MR1970804}}]
\label{prop:hilbert-function-of-WLP-SLP}
Let $h = (h_0, h_1, \ldots, h_c)$ be a sequence of positive integers.
The following three conditions are equivalent.
\begin{enumerate}
  \item[(i)]
  $h$ is a  Hilbert function of some graded algebra with the SLP,
  \item[(ii)]
  $h$ is a  Hilbert function of some graded algebra with the WLP,
  \item[(iii)]
  $h$ is a unimodal O-sequence,
  and the sequence $\Delta h$ is an O-sequence.
\end{enumerate}
\qed
\end{proposition}
We introduce the notion of $x_n$-chains.
Let $R = \PPn$ be the polynomial ring over a field $K$.
We take any term order $\sigma$ on $R$.
When a graded ideal $I$ of $R$ is given,
a monomial which does not belong 
to the initial ideal $\In_\sigma(I)$ of $I$
is called a {\it standard monomial} with respect to $I$.
The standard monomials with respect to $I$ is the same as
the standard monomials with respect to $\In_\sigma(I)$.
The standard monomials span the quotient ring $R/I$
as a $K$-vector space
(\cite[Theorem 4.2.3]{MR1251956}).
\begin{definition and remark}
\label{defn:xn-chain} 
Let $I \subset R$ be a graded Artinian ideal.
An {\it $x_n$-chain} with respect to $I$ and a term order $\sigma$
is defined as the sequence 
\begin{align*}
  \{ u x_n^s \;|\; s \ge 0,\, u x_n^s \not\in \In_\sigma(I) \}
  \qquad
  (u \in R,\; x_n \nmid u)
\end{align*}
consisting of standard monomials,
where $u \in R$ is a monomial not divisible by $x_n$.
Note that the $x_n$-chains with respect to $I$ is the same as
that with respect to $\In_\sigma(I)$. 

The set of the standard monomials decomposes into
disjoint $x_n$-chains,
and we call this decomposition the {\it $x_n$-chain decomposition}
of the standard monomials with respect to $I$ and $\sigma$.
Note that the $x_n$-chain decomposition 
with respect to $I$ and $\sigma$
is uniquely determined by $I$ and $\sigma$,
and that the  $x_n$-chain decomposition 
with respect to $I$ and $\sigma$
is the same as that with respect to $\In_\sigma(I)$.
Moreover, if $I$ is a monomial ideal, 
then the $x_n$-chain decomposition does not depend on term orders.
\end{definition and remark}
For example, consider 
$I = (x_1^2, x_1x_2, x_2^2, x_2x_3^2) + (x_1,x_2,x_3)^4$
in $R = \PPiii$.
The $x_3$-chain decomposition of the standard monomials is as follows:
\begin{align}
  \label{eq:SLP-std-mono}
  \begin{array}{cccccc}
    {} & x_1, & x_1x_3, & x_1x_3^2, \\
    {} & x_2, & x_2x_3, \\
    1, & x_3, & x_3^2, & x_3^3.
  \end{array}
  \qquad
  \text{(each row forms an $x_3$-chain)}
\end{align}
For 
$I = (x_1^2, x_1x_2, x_2^3, x_1x_3^2, x_2^2x_3^2, x_2x_3^3, x_3^5)$, 
the $x_3$-chain decomposition of the standard monomials is as follows:
\begin{align}
  \label{eq:non-SLP-std-mono}
  \begin{array}{ccccc}
    {}  &  {}    &  x_2^2,   &  x_2^2x_3, \\
    {}  &  x_1,  &  x_1x_3,  \\
    {}  &  x_2,  &  x_2x_3,  &  x_2x_3^2, \\
    1,  &  x_3,  &  x_3^2,   &  x_3^3,     &  x_3^4. 
  \end{array}
  \qquad
  \text{(each row forms an $x_3$-chain)}
\end{align}

We can also consider $x_i$-chains for $i < n$,
but we do not need them in this paper.

Here we give a necessary and sufficient condition
for a quotient ring $R/I$ to have the SLP or the WLP
in terms of $x_n$-chains
introduced in Definition~\ref{defn:xn-chain}.
This condition is only a paraphrase of the ordinary conditions
in terms of minimal generators and Hilbert functions.
Our condition, however,
concretely describes the structure of standard monomials 
in contrast to the ordinary conditions.
In particular,
an advantage of our conditions appears 
when we construct a monomial ideal
for which the quotient ring has the SLP and 
has a given Hilbert function.
\begin{definition}
\label{defn:Lefschetz-condition} 
Let $R = \PPn$, and $I \subset R$ a graded Artinian ideal,
and fix any term order on $R$.

(i)
For the $x_n$-chain decomposition of the standard monomials,
we call the following condition the {\it strong Lefschetz condition}
(the {\it SL condition}, for short):
\begin{quote}
  For any two $x_n$-chains
  $u, u x_n, u x_n^2, \ldots, u x_n^s$ and
  $v, v x_n, v x_n^2, \ldots, v x_n^t$,
  where $u, v \in R$ are monomials not divisible by $x_n$,
  if $\deg u < \deg v$, then $\deg u x_n^s \ge \deg v x_n^t$.
\end{quote}
Roughly speaking,
this condition means that
a chain which begins at a lower degree 
continues until a higher degree.

(ii)
For the $x_n$-chain decomposition of the standard monomials,
we call the following condition the {\it weak Lefschetz condition}
(the {\it WL condition}, for short):
\begin{quote}
Let $u_2$ be the minimum integer satisfying $h_{u_2-1} > h_{u_2}$
(see Equation~(\ref{eq:def-of-ui})),
where $(h_0, h_1, \ldots, h_c)$ denotes the Hilbert function of $R/I$.
Then every $x_n$-chain starts 
at a degree less than or equal to $u_1$,
and ends in a degree greater than or equal to $u_2-1$.
\end{quote}
In other words,
if an $x_n$-chain ends in degree $d$,
then no $x_n$-chain starts from degrees
greater than $d$.
\end{definition}
For example,
the $x_3$-chain decomposition in Figure~(\ref{eq:SLP-std-mono})
satisfies the SL condition.
The $x_3$-chain decomposition in Figure~(\ref{eq:non-SLP-std-mono})
does not satisfy the SL condition, 
but satisfies the WL condition.

The following lemma relates the Lefschetz properties with
the Lefschetz conditions for the $x_n$-chain decomposition.
This lemma is fundamental
in studying the standard monomials 
whose quotient rings have the SLP or the WLP.
\begin{lemma}
\label{lem:lefschetz-condition} 
Let $I$ be a graded Artinian ideal of $R=\PPn$,
and consider the graded reverse lexicographic order on $R$.
The following conditions are equivalent:
\begin{itemize}
  \item[(i)]
  $(R/I, x_n)$ has the SLP (resp. WLP).
  \item[(ii)]
  $(R/\inrevlex(I), x_n)$ has the SLP (resp. WLP).
  \item[(iii)]
  The $x_n$-chain decomposition of the standard monomials
  with respect to $I$ satisfies the SL condition
  (resp. WL condition). 
\end{itemize}
\end{lemma}
\begin{proof}
It is obvious that (ii) and (iii) are equivalent.
We show the equivalence of (i) and (ii).
Let $\In(I)$ denote the initial ideal of $I$
with respect to the graded reverse lexicographic order.
We have $ \In(I:x_n^s) = \In(I):x_n^s \quad (s>0)$
\cite[Proposition 15.12]{MR1322960}.
By the following linear isomorphisms
\begin{gather*}
\Ker( R/I \overset{\times x_n^s}{\to} R/I ) = I:x_n^s / I,
\\
\Ker( R/\In(I) \overset{\times x_n^s}{\to} R/\In(I) ) = 
\In(I):x_n^s / \In(I),
\end{gather*}
we have
\begin{align*}
  \rank( \times x_n^s: (R/I)_i \to (R/I)_{i+s} )
  =
  \rank( \times x_n^s: (R/\In(I))_i \to (R/\In(I))_{i+s} )
  \\
  (i \ge 0, \; s > 0),
\end{align*}
where `rank' means the rank of a linear map,
and the $i$-th graded components are denoted as $(R/I)_i$ and so on.
For a graded ideal $J$ of $R$,
$(R/J, x_n)$ has the SLP (resp. WLP),
if and only if the linear map
$\times x_n^s: (R/J)_i \to (R/J)_{i+s}$ is full-rank
for every $i\ge0$ and $s>0$ (resp.~$s=1$).
Therefore it follows from the formula above
that (i) and (ii) of the lemma are equivalent.
\end{proof}

\subsection{Almost revlex ideals and the SLP}
\label{subsec:almost-revlex}

\begin{definition}
\label{defn:Borel-fixed-strongly-stable}
Let $I$ be an ideal of $R = \PPn$.
For an invertible matrix $g = (g_{i j})_{1\le i,j \le n} \in GL(n;K)$,
the transform $g(I)$ of the ideal $I$ is defined by the image of the action
$g x_j = \sum_{i=1}^n g_{i j} x_i$.
An ideal $I$ of $R$ is said to be {\it Borel-fixed},
if $g(I) = I$ for every upper triangular matrix $g \in GL(n;K)$.

A monomial ideal $I$ is said to be
{\it strongly stable} (resp. {\it stable})
if $I$ satisfies the following condition:
\begin{gather}
  \label{eq:condition-of-Borel}
  \parbox{0.8\textwidth}{
    for each monomial $u \in I$ and an index $j$
    (resp. the maximum index $j$) satisfying $x_j | u$,
    the monomial $x_i u/x_j$ belongs to $I$ for every $i < j$.}
\end{gather}

It is known that strongly stable ideals are Borel-fixed
for any characteristic,
and that Borel-fixed ideals are strongly stable
when the characteristic is zero
(see \cite[15.9]{MR1322960}, e.g.).
\end{definition}

We recall a result of Wiebe \cite{MR2111103}.
This lemma shows that the $x_n$-chain decomposition is a useful tool 
in studying stable ideals 
whose quotient rings have the SLP or the WLP.
\begin{lemma}[Wiebe, {\cite[Lemma 2.7]{MR2111103}}]
\label{lem:wiebe-2.7}
If $I$ is an Artinian stable ideal of $R=\PPn$,
then the following two conditions are equivalent:
\begin{itemize}
  \item[(i)]
  $R/I$ has the SLP (resp. WLP),
  \item[(ii)]
  $x_n$ is a strong (resp. weak) Lefschetz element on $R/I$.
\end{itemize}
\qed
\end{lemma}
We show another lemma, 
which also relates stable ideals with $x_n$-chains.
This lemma is a special case of 
Lemma~\ref{lem:m-full-is-xn-chain} proved below,
and we omit the proof.
\begin{lemma}
\label{lem:stable-is-xn-chain}
Let $I$ be an Artinian stable ideal of $\PPn$ ,
and $u, u x_n, \ldots, u x_n^{s-1}$ an $x_n$-chain
with respect to $I$,
where $u$ is a monomial not divisible by $x_n$, and $s \ge 1$.
Then $u x_n^s$ is a member of the minimal generators of $I$.
\qed
\end{lemma}
We define the notion of almost revlex ideals
(almost revlex-segment ideals).
Let $R = \PPn$ be the polynomial ring over a field 
of characteristic zero.
Let $>_{\text{revlex}}$ denote the graded reverse lexicographic order.
\begin{definition}
\label{defn:almost-revlex-ideal} 
(i) A monomial ideal $I$ is called a {\it revlex ideal},
if the following condition holds:
\begin{quote}
  for each monomial $u \in I$,
  every monomial $v$ with $\deg v = \deg u$ and
  $v >_{\text{revlex}} u$ belongs to $I$.
\end{quote}

(ii) A monomial ideal $I$ is called an {\it almost revlex ideal},
if the following condition holds:
\begin{quote}
  for each monomial $u$ in the minimal generating set of $I$,
  every monomial $v$ with $\deg v = \deg u$ 
  and $v >_{\text{revlex}} u$ belongs to $I$.
\end{quote}
\end{definition}
\begin{remark}
\label{rem:almost-revlex-is-stable}
First it is clear that 
\begin{itemize} 
  \item[(i)] revlex ideals are almost revlex ideals. 
\end{itemize}
Second,
\begin{itemize} 
  \item[(ii)] if two almost revlex ideals have 
  the same Hilbert function, then they are equal,
\end{itemize}
since one can determine the minimal generators from low degrees
using a given Hilbert function.
\end{remark}
In addition to Remark~\ref{rem:almost-revlex-is-stable},
we have the following lemma.
\begin{lemma}
\label{lem:almost-revlex-is-Borel}
Almost revlex ideals are Borel-fixed.
\end{lemma}
\begin{proof}
Let $w$ be a monomial of an almost revlex ideal $I$,
and suppose that $x_j \mid w$.
We show $(x_i w)/x_j \in I $ for all $i < j$.
Take a monomial $u$ in the minimal generating set of $I$
such that $w = v u$ for some monomial $v$.
If $x_j\mid v$ then $(x_i w)/x_j = ((x_i v)/x_j)u \in I$.
If $x_j \mid u$ then $(x_i u)/x_j \in I$, 
since $(x_i u)/x_j >_{revlex} u$ 
and $I$ is an almost revlex ideal.
Hence $(x_i w)/x_j = v((x_i u)/x_j) \in I$.
\end{proof}
The notion `almost revlex' is paraphrased
in terms of $x_n$-chains as follows.
\begin{proposition}
\label{prop:paraphrase-almost-revlex} 
Let $I$ be a monomial Artinian ideal of $R = \PPn$,
and set $\Rbar = \PPnm$.
Then $I$ is an almost revlex ideal, if and only if
$I$ satisfies the following two conditions:
\begin{itemize}
  \item[(i)]
  $I \cap \Rbar$ is almost revlex,
  \item[(ii)]
  For two standard monomials $u <_{\text{revlex}} v$
  not divisible by $x_n$,
  the ending degree of the $x_n$-chain starting with $u$
  is greater than or equal to 
  that of the $x_n$-chain starting with $v$.
\end{itemize}
In the condition (ii), 
revlex denotes the graded reverse lexicographic order,
and note that the degrees of $u$ and $v$ are not necessarily equal.
\end{proposition}
Since the condition (ii)
of Proposition~\ref{prop:paraphrase-almost-revlex}
in the case where  $\deg u < \deg v$ is nothing but the SL condition,
we have the following corollary.
\begin{corollary}
\label{cor:almost-revlex-has-SLP}
Let $I \subset R = \PPn$ be an Artinian almost revlex ideal.
Then $(R/I, x_n)$ has the SLP.
\qed
\end{corollary}
\begin{proof}%
[{Proof of Proposition~\ref{prop:paraphrase-almost-revlex}}]
Let $I$ be an almost revlex ideal.
It is clear that $I \cap \Rbar$ is almost revlex, and we prove (ii).
Let $u$ and $v$ be standard monomials not divisible by $x_n$
with $u <_{\text{revlex}} v$.
Note that $\deg u \le \deg v$ in this case.
Suppose that $v x_n^{s-\deg v}$ is standard ($s \ge \deg v$).
It suffices to show that $u x_n^{s-\deg u}$ is standard
for proving (ii).
Assume that $u x_n^{s-\deg u} \in I$, 
then there exists an integer $t$ ($\deg u < t \le s$) 
such that $u x_n^{t-\deg u}$ is a minimal generator
thanks to Lemma~\ref{lem:stable-is-xn-chain}.
Since $I$ is almost revlex,
if $t \ge \deg v$ then $v x_n^{t-\deg v} \in I$,
and hence $v x_n^{s-\deg v} \in I$.
This is a contradiction, and therefore $t < \deg v$.
In this case, all the monomials in $x_1, \ldots, x_{n-1}$
of degree $t$ belong to $I$,
since $ux_n^{t-\deg(u)}$ is a minimal generator 
and $I$ is almost revlex.
This yields that $v \in I$, 
and this is also a contradiction.
We therefore have $u x_n^{s-\deg u} \not\in I$, and proved (ii).

Conversely, let $I$ be a monomial ideal of $R$
satisfying (i) and (ii).
We prove that $I$ is almost revlex.
Let $u x_n^s$ ($s \ge 0$) be a minimal generator of $I$,
where $u$ is a monomial not divisible by $x_n$.
We have to show that every monomial of the same degree as $u x_n^s$
which is greater than $u x_n^s$ with respect to $<_{\text{revlex}}$
belongs to $I$.
We first consider the case of $s=0$.
Monomials of the same degree as $u$ which is greater than $u$
are monomials in $x_1, x_2, \ldots, x_{n-1}$.
Such monomials are in $I \cap \Rbar$,
since $I \cap \Rbar$ is almost revlex.
Hence such monomials belong to $I$.
Second we consider the case of $s > 0$.
In this case $u$ is a standard monomial.
Take a monomial $v x_n^t$ ($t \ge 0$) of the same degree as $u x_n^s$,
where $v$ is not divisible by $x_n$,
and $v x_n^t >_{\text{revlex}} u x_n^s$.
We have to show that $v x_n^t \in I$.
The relation $v x_n^t >_{\text{revlex}} u x_n^s$ implies 
that $v >_{\text{revlex}} u$ ($t=s$)
or $\deg v > \deg u$ ($t<s$).
In both cases, we have $v >_{\text{revlex}} u$.
Then $v x_n^t$ can not be a standard monomial 
because of the condition (ii).
Hence $v x_n^t \in I$.
We thus have proved that $I$ is an almost revlex ideal.
\end{proof}

\subsection{Uniqueness of Borel-fixed ideals and 
generic initial ideals in $\PPiii$}

For a given O-sequence, it is known that
a Borel-fixed ideal of $\PPii$,
whose quotient ring has the O-sequence as the Hilbert function,
is unique.
It is the unique lex-segment ideal determined 
by the Hilbert function.
Moreover we have the following theorem,
which gives the uniqueness of Borel-fixed ideals
for $n = 3$, where the quotient rings have the SLP.
\begin{theorem}
\label{thm:3var-uniqueness-of-Borel-SLP} 
Let $R = \PPiii$ be the polynomial ring 
over a field of characteristic zero,
and $I$ an Artinian Borel-fixed ideal of $R$ where $R/I$ has the SLP.
Then the ideal $I$ is the unique almost revlex ideal 
for the Hilbert function.
In particular,
the ideal $I$ is uniquely determined only by the Hilbert function.
\end{theorem}
\begin{proof}
We check the conditions (i) and (ii) of 
Proposition~\ref{prop:paraphrase-almost-revlex}.
First, $I \cap \PPii$ is Borel-fixed, and hence almost revlex,
since Borel-fixed ideals are revlex ideals in $\PPii$.
Thus the condition (i) of 
Proposition~\ref{prop:paraphrase-almost-revlex} holds.

Second we check the condition (ii).
For two standard monomials $u <_{\text{revlex}} v$
not divisible by $x_3$,
whose degrees are not necessarily equal,
we prove that
the ending degree of the $x_3$-chain starting with $u$ 
is greater than or equal to that of the $x_3$-chain starting with $v$.
In the case where $\deg u < \deg v$,
the condition (ii) of 
Proposition~\ref{prop:paraphrase-almost-revlex} is 
nothing but the SL condition as mentioned 
before Corollary~\ref{cor:almost-revlex-has-SLP}.
Since $(R/I, x_3)$ has the SLP, this condition holds. 
Next suppose that $\deg u = \deg v = k$,
and let $u = x_1^a x_2^{k-a}$, $v = x_1^b x_2^{k-b}$, and $a < b$.
If $v x_3^s$ is standard then $u x_3^s$ is also standard,
since $I$ is Borel-fixed.
This means again that the condition (ii) of 
Proposition~\ref{prop:paraphrase-almost-revlex} holds.

Thus the conditions (i) and (ii) of 
Proposition~\ref{prop:paraphrase-almost-revlex} hold,
and hence $I$ is almost revlex.
Finally $I$ is determined only by the Hilbert function,
by Remark~\ref{rem:almost-revlex-is-stable} (ii).
\end{proof}
If the number of variables is more than three,
then Theorem~\ref{thm:3var-uniqueness-of-Borel-SLP} does not hold.
In the middle of the proof
we take two standard monomials  $u <_{\text{revlex}} v$ 
of the same degree not divisible by $x_3$,
and concluded that $u x_3^s$ is standard if $v x_3^s$ is standard.
But this argument fails when the number of the variables is
greater than three.
Take $u <_{\text{revlex}} v$ as
$u = x_1 x_3$ and $v = x_2^2$ in $\PPiv$, for example.
In this case, we can not conclude 
from the condition `Borel-fixed'
that $u x_4^s$ is standard even if $v x_4^s$ is standard.
See Example~\ref{ex:4var-non-symmetric-Borel-2SLP-is-not-unique}
for a counterexample in the case of four variables.

\begin{example}
\label{ex:3var-Borel-WLP-is-not-unique}
For a given Hilbert function, Borel-fixed ideals in $R = \PPiii$,
where the quotient rings have the WLP,
are not unique in general.
Consider the following distinct Borel-fixed ideals in $R$:
\begin{align*}
  I &= (x_1^2, x_1x_2, x_2^3, x_2^2x_3, x_1x_3^3, x_2x_3^3, x_3^5),
  \\
  J &= (x_1^2, x_1x_2, x_2^3, x_1x_3^2, x_2^2x_3^2, x_2x_3^3, x_3^5).
\end{align*}
We can easily check that 
both $(R/I, x_3)$ and $(R/J, x_3)$ have the WLP
and the same symmetric Hilbert function
$h = (1,3,4,3,1)$.
Note that $R/I$ has the SLP, but $R/J$ does not.
The standard monomials with respect to $J$ is given in
Figure~(\ref{eq:non-SLP-std-mono}).
\end{example}
\begin{example}
\label{ex:3var-stable-SLP-is-not-unique}
For given Hilbert functions,
stable ideals in $R = \PPiii$,
where the quotient rings have the SLP,
are not unique in general.
Consider the following distinct stable ideals in $R$:
\begin{align*}
  I &= (x_1^2, x_1x_2, x_2^2, x_1x_3^2) + (x_1,x_2,x_3)^4,
  \\
  J &= (x_1^2, x_1x_2, x_2^2, x_2x_3^2) + (x_1,x_2,x_3)^4.
\end{align*}
We can easily check that 
both $(R/I, x_3)$ and $(R/J, x_3)$ have the SLP
and the same Hilbert function $h = (1,3,3,2)$.
Note that $I$ is Borel-fixed, but $J$ is not.
The standard monomials with respect to $J$ is given in
Figure~(\ref{eq:SLP-std-mono}).
\end{example}

The following is an immediate corollary to
Proposition~\ref{prop:hilbert-function-of-WLP-SLP},
Corollary~\ref{cor:almost-revlex-has-SLP} ,
Theorem~\ref{thm:3var-uniqueness-of-Borel-SLP}
and Lemma~\ref{lem:wiebe-2.8} below.
\begin{corollary}
\label{cor:hilbert-fnct-of-3SLP}
Let $R = \PPiii$ 
and $h = (1, 3, h_2, h_3, \ldots, h_c)$ an O-sequence.
The following three conditions are equivalent:
\begin{itemize}
  \item[(i)]
  $h$ is a Hilbert function of $R/I$ 
  for some almost revlex ideal $I$ of $R$,
  \item[(ii)]
  $h$ is a Hilbert function of some graded algebra with the SLP,
  \item[(iii)]
  $h$ is a Hilbert function of some graded algebra with the WLP,
  \item[(iv)]
  $h$ is a unimodal O-sequence, 
  and $\Delta h$ is an O-sequence. 
  \qed
\end{itemize}
\end{corollary}
In the rest of this section,
we study generic initial ideals in $\PPiii$.
We recall the definition of generic initial ideals.
Fix any term order $\sigma$ on
the polynomial ring $R = K[x_1, x_2, \ldots, x_n]$
over a field of characteristic zero.
For a graded ideal $I$ of $R$,
there exists a Zariski open subset $U \subset GL(n; K)$
such that the initial ideals of $g(I)$ are equal to each other
for any $g \in U$.
This initial ideal is uniquely determined,
called the {\it generic initial ideal} of $I$,
and denoted by $\gin_\sigma(I)$.
It is known that generic initial ideals 
are Borel-fixed with respect to any term order
(see \cite[15.9]{MR1322960}, e.g.).

Thus we have results on generic initial ideals of ideals
whose quotient rings have the SLP,
as an easy consequence 
of Theorem~\ref{thm:3var-uniqueness-of-Borel-SLP}.
We first recall another result of Wiebe.
We write simply by $\gin(I)$
the generic initial ideal of $I$ 
with respect to the graded reverse lexicographic order from now on.
\begin{lemma}[Wiebe, {\cite[Proposition~2.8]{MR2111103}}]
\label{lem:wiebe-2.8}
Take the graded reverse lexicographic order on 
the polynomial ring $R = \PPn$
over a field $K$ of characteristic zero.
Let $I$ be a graded Artinian ideal of $R$.
Then $R/I$ has the SLP if and only if $R/\gin(I)$ has the SLP.
\qed
\end{lemma}
We thus have the following theorem by
Lemma~\ref{lem:wiebe-2.7},
Theorem~\ref{thm:3var-uniqueness-of-Borel-SLP} and
Lemma~\ref{lem:wiebe-2.8}.
\begin{theorem}
\label{thm:3var-gin} 
Let $R = \PPiii$ be the polynomial ring 
over a field of characteristic zero,
and consider the graded reverse lexicographic order on $R$.
Let $I$ be a graded Artinian ideal of $R$,
and suppose that $R/I$ has the SLP.
Then the generic initial ideal $\gin(I)$ is 
the unique almost revlex ideal 
for the same Hilbert function as $R/I$.
In particular, 
$\gin(I)$ is uniquely determined by the Hilbert function.
\qed
\end{theorem}
Cimpoea\c{s} \cite{math-0610649} shows that
the generic initial ideals of complete intersections of height three
whose quotient rings have the SLP
are almost revlex ideals.
Theorem~\ref{thm:3var-gin} is an improvement of this result.
In addition, Theorem~\ref{thm:3var-gin} is an improvement 
of a result of Ahn, Cho and Park \cite{math-0610733}.
They prove that
the generic initial ideals of ideals in $\PPiii$
whose quotient rings have the SLP
are determined by their graded Betti numbers. 
\begin{remark}
\label{rem:3var-gin-of-WLP-is-not-unique}
For a graded ideal $I \subset \PPiii$, 
if its quotient ring has the WLP and does not have the SLP,
then its generic initial ideal is not determined
by the Hilbert function in general.
A counterexample is already given in 
Example~\ref{ex:3var-Borel-WLP-is-not-unique},
since the generic initial ideal of a Borel-fixed ideal $I$
is equal to $I$ itself.
\end{remark}
\begin{example}
\label{ex:3var-ci-SLP}
We give four examples of complete intersection in $R = \PPiii$
whose quotient rings have the SLP.
The generic initial ideals of these ideals are
the unique almost revlex ideals with corresponding 
Hilbert functions by Theorem~\ref{thm:3var-gin}.

(i) Let $I = (f, g, \ell^r) \subset R$,
where $f$ and $g$ are any homogeneous polynomials of $R$,
and $\ell$ is any homogeneous polynomial of degree one.
In this case, if $I$ is a complete intersection, then
$R/I$ has the SLP
(\cite[Example 6.2]{harima-watanabe-jpaa}).

(ii) Let $e_1$, $e_2$ and $e_3$ be the elementary symmetric
functions in three variables, where $\deg(e_i) = i$.
Let $r$ and $s$ be positive integers, where $r$ divides $s$.
Then, the quotient ring of the ideal
$I=(e_1(x_1^r, x_2^r, x_3^r),
  e_2(x_1^r, x_2^r, x_3^r),
  e_3(x_1^s, x_2^s, x_3^s))$
of $R$ has the SLP
(\cite[Example 6.4]{harima-watanabe-jpaa}).

(iii) Let $p_i$ be the power sum symmetric function of degree $i$
in three variables,
and $a$ be a positive integer.
Then, the quotient ring of the ideal
$I =  (p_a, p_{a+1}, p_{a+2})$ of $R$ has the SLP
(\cite[Proposition 7.1]{harima-watanabe-ja}).

(vi) Let $I = (e_2, e_3, f) \subset R$,
where $f$ is any homogeneous polynomial of $R$.
In this case, if $I$ is a complete intersection, then
$R/I$ has the SLP
(\cite[Proposition 3.1]{harima-watanabe-ja}).
\end{example}

\section{Generic initial ideals in $\PPn$ and the $k$-SLP}

Suppose that a graded Artinian algebra $A$ has the SLP (resp. WLP),
and $\ell \in A$ is a Lefschetz element.
If the graded algebra $A/(\ell)$ again has the SLP (resp. WLP),
then we say that $A$ has the 2-SLP (resp. 2-WLP).
We define the notion of the $k$-SLP and the $k$-WLP
recursively (Definition~\ref{defn:k-SLP-k-WLP}).
A characterization of the Hilbert functions 
of graded Artinian algebras having the $k$-SLP or the $k$-WLP
is given in Proposition~\ref{prop:hilbert-fnct-of-nSLP}.
Moreover, the Hilbert functions of quotient rings
by almost revlex ideals are determined in terms of the $n$-SLP
in Proposition~\ref{prop:hilbert-fnct-of-almost}.

The main goal of this section is
Theorem~\ref{thm:gin-of-symmetric-nSLP}:
Let $I \subset \PPn$ be a graded Artinian ideal
whose quotient ring has the $n$-SLP,
and every $k$-th difference of the Hilbert function 
is quasi-symmetric.
The generic initial ideal of $I$
with respect to the graded reverse lexicographic order
is the unique almost revlex ideal 
for the same Hilbert function as $\PPn/I$.

\subsection{$k$-SLP and $k$-WLP}
\label{subsec:k-SLP-k-WLP}
The first author heard from J. Watanabe that
the following notion, the `$k$-SLP' and the `$k$-WLP', 
has been introduced by A.~Iarrobino 
in a private conversation with J. Watanabe in 1995. 
\begin{definition}
\label{defn:k-SLP-k-WLP}
Let $A = \bigoplus_{i=0}^c A_i$ be a graded Artinian $K$-algebra,
and $k$ a positive integer.
We say that $A$ has the {\it $k$-SLP} (resp. {\it $k$-WLP})
if there exist linear elements $g_1, g_2, \ldots, g_k \in A_1$
satisfying the following two conditions.
\begin{itemize}
  \item[(i)]
  $(A, g_1)$ has the SLP (resp. WLP),
  \item[(ii)]
  $(A/(g_1, \ldots, g_{i-1}), g_i)$ has the SLP (resp. WLP)
  for all $i = 2, 3, \ldots, k$.
\end{itemize}
In this case, we say that 
$(A, g_1, \ldots, g_k)$ has the $k$-SLP (resp. $k$-WLP).
Note that a graded algebra with the $k$-SLP (resp. $k$-WLP)
has the $(k-1)$-SLP (resp. $(k-1)$-WLP).
\end{definition}
\begin{remark}
\label{rem:k-SLP-k-WLP}
From Theorem~4.4 in \cite{MR1970804}, 
one knows that 
all graded $K$-algebras $K[x_1]/I$ and $\PPii/I$ have the SLP. 
Hence the following conditions are equivalent 
for a graded Artinian algebra $A=\PPn/I$,
where $I \subset (x_1, \ldots, x_n)^2$.
\begin{itemize}
  \item[(i)]
  $A$ has the $n$-WLP (resp. the $n$-SLP),
  \item[(ii)]
  $A$ has the $(n-1)$-WLP (resp. the $(n-1)$-SLP),
  \item[(iii)]
  $A$ has the $(n-2)$-WLP (resp. the $(n-2)$-SLP). 
\end{itemize}
In particular, graded algebras $\PPiii/I$ with the SLP (resp. WLP)
has the 3-SLP (resp. 3-WLP) automatically.
\end{remark}
\begin{example}
\label{ex:almost-revlex-is-nSLP}
For every Artinian almost revlex ideal $I$ of $R=\PPn$
where $I \subset (x_1,\ldots,x_n)^2$,
the quotient ring $R/I$ has the $n$-SLP. 
In particular, for every revlex ideal $I$ of $R$
where $I\subset (x_1,\ldots,x_n)^2$, 
the quotient ring $R/I$ has the $n$-SLP. 
\end{example}
\begin{proof}
Let $\Rbar = \PPnm$.
If $I$ is an almost revlex ideal of $R$, 
then it follows from Corollary~\ref{cor:almost-revlex-has-SLP}
that $(R/I, x_n)$ has the 1-SLP,
and $I \cap \Rbar$ is again an almost revlex ideal.
This shows that
$I+(x_n)/(x_n)$ is an almost revlex ideal
as an ideal of $R/(x_n) \simeq \Rbar$,
and hence $(R/I, x_n, x_{n-1})$ has the 2-SLP.
Repeating this argument, we obtain that 
$(R/I, x_n, x_{n-1}, \ldots, x_1)$  has the $n$-SLP.
\end{proof}
Example~\ref{ex:almost-revlex-is-nSLP} shows that
the class of Hilbert functions for almost revlex ideals
is a subset of that for ideals with the $n$-SLP.
In fact, 
Proposition~\ref{prop:hilbert-fnct-of-almost}
shows that these two classes coincide.
We also determine the class of Hilbert functions
for ideals with the $k$-SLP 
in Proposition~\ref{prop:hilbert-fnct-of-nSLP}.
\par
Let $\gin(I)$ denote the generic initial ideal of $I$
with respect to the graded reverse lexicographic order.
The following proposition is an analogue of 
Wiebe's result (Lemma~\ref{lem:wiebe-2.8})
\cite[Proposition~2.8]{MR2111103}.
\begin{proposition}
\label{prop:I-is-kSLP-iff-Gin(I)-is-kSLP}
Let $I$ be a graded Artinian ideal of $R=K[x_1,\ldots,x_n]$,
and let $1\leq k\leq n$. 
The following two conditions are equivalent: 
\begin{itemize}
  \item[(i)]
  $R/I$ has the $k$-WLP (resp. the $k$-SLP),
  \item[(ii)]
  $(R/\gin(I), x_n, x_{n-1}, \ldots, x_{n-k+1})$
  has the $k$-WLP (resp. the $k$-SLP). 
\end{itemize}
\end{proposition}
\begin{proof}
We can give a proof using a similar idea to the proof of 
Proposition~2.8 in \cite{MR2111103}. 
First we show that (i) and (ii) are equivalent for the $k$-WLP. 
Let $1\leq i\leq k$. 
Lemma~1.2 in \cite{MR1948090} says that 
the Hilbert function of $R/(I + (g_1,\ldots,g_i))$ 
for generic linear forms $g_1,\ldots,g_i$ 
is equal to the Hilbert function of
$R/(\gin(I) + (x_n, x_{n-1}, \ldots x_{n-i+1}))$.
Noting the last variable is a Lefschetz element 
of quotient rings by Borel-fixed ideals
when they have the WLP,
it follows from Remark~\ref{rem:WLP-and-SLP} (i) 
that this means the equivalence of (i) and (ii). 

Next we show that (i) and (ii) are equivalent for the $k$-SLP. 
Let $1\leq i\leq k$. 
By slightly generalizing the proofs of Lemma~1.2 in \cite{MR1948090}
and Proposition~2.8 in \cite{MR2111103}, 
it follows that 
the Hilbert function of $R/(I + (g_1,\ldots,g_{i-1},g_i^s))$
is equal to the Hilbert function of 
$R/(\gin(I)+(x_n, \ldots, x_{n-i+2}, x_{n-i+1}^s))$ 
for generic linear forms $g_1,\ldots,g_i$ and $s \ge 1$. 
Noting that the last variable is a Lefschetz element 
of quotient rings by Borel-fixed ideals when they have the SLP, 
it follows from Remark~\ref{rem:WLP-and-SLP} (ii) 
that this means the equivalence of (i) and (ii) for the $k$-SLP. 
\end{proof}
\begin{remark}
\label{rem:generic-elements-of-Borel-kSLP}
Let $I$ be a graded Artinian ideal of $R = \PPn$
whose quotient ring $R/I$ has the $k$-SLP (resp. $k$-WLP).
If $I$ is Borel-fixed, 
then $\gin(I) = I$.
Hence, by Proposition~\ref{prop:I-is-kSLP-iff-Gin(I)-is-kSLP},
$(R/I, x_n, x_{n-1}, \ldots, x_{n-k+1})$ has 
the $k$-SLP (resp. $k$-WLP).
In other words, we can choose the variable $x_{n-i+1}$
as a Lefschetz element on $R/(I, x_n, x_{n-1}, \ldots, x_{n-i+2})$
for all $i = 1,2,\ldots, k$.
\end{remark}

\subsection{Correspondence of Borel-fixed ideals and the $k$-SLP}

We consider graded Artinian ideals $I \subset R = \PPn$ 
such that $(R/I, x_n)$ has the SLP,
and define the map $\varphi(I) = I \cap \PPnm$.
We study behaviors of this map when restricting 
to the Borel-fixed ideals or the almost revlex ideals.

For an O-sequence $h$, let $\cM_n(h)$ denote 
the set of the monomial ideals of $R$ 
for which $R/I$ has the Hilbert function $h$.
We need the following definition 
for Lemma~\ref{lem:intersection-with-PPnm}.
\begin{definition}
\label{defn:quasi-symmetric}
A unimodal sequence $h = (h_0, h_1, \ldots, h_c)$ of positive integers
is said to be {\it quasi-symmetric}, 
if the following condition holds:
\begin{quote}
  Let $h_i$ be the maximum of $\{h_0, h_1, \ldots, h_c \}$.
  Then every integer $h_j$ ($j > i$) is equal to one of 
  $\{ h_0, h_1, \ldots, h_i \}$.
\end{quote}
In particular,
unimodal symmetric sequences are quasi-symmetric.
\end{definition}
\begin{lemma}
\label{lem:intersection-with-PPnm}
Set $R = \PPn$ and $\Rbar = \PPnm$.
Let $h = (1, n, h_2, \ldots, h_c)$ be an O-sequence
which can be a Hilbert function of a graded algebra with the SLP.
Consider the mapping
\begin{gather*}
  \begin{array}{cccc}
    \varphi: &
    \{ I \in \cM_n(h) \;;\; \text{$(R/I,x_n)$ has the SLP} \}
    & \to &
    \cM_{n-1}(\Delta h)
    \\ \\
    & I &\mapsto& I \cap \Rbar.
  \end{array}
\end{gather*}
Note that the domain of the map $\varphi$ is not empty
thanks to Lemmas~\ref{lem:wiebe-2.7} and \ref{lem:wiebe-2.8}.
We have the following.
\begin{itemize}
  \item[(i)]
  $\varphi$ is surjective,
  \item[(ii)]
  $\varphi$ is bijective, if $h$ is quasi-symmetric.
\end{itemize}
\end{lemma}
\begin{proof}
We first note that the map $\varphi$ is well-defined
thanks to Remark~\ref{rem:WLP-and-SLP}.

(i)
For a graded Artinian ideal $\Ibar \in \cM_{n-1}(\Delta h)$,
we construct $I \in \cM_n(h)$ for which $(R/I, x_n)$ has the SLP.
For this purpose, we construct the set $M$  of the standard monomials
with respect to $I$,
since the monomials not belonging to $M$ generate $I$.
Let $v_1, v_2, \ldots, v_m \in \Rbar$ be the standard monomials 
with respect to $\Ibar$,
where $\deg v_i \le \deg v_{i+1}$.
In each degree, the ordering of $v_i$'s is arbitrary.
Let $p_1, p_2, \ldots, p_m$ be the lengths of the $x_n$-chains
with respect to $I$,
where duplicated lengths repeatedly appear,
and $p_i \ge p_{i+1}$.
Note that $p_i$'s are determined only by the O-sequence $h$,
since $R/I$ should have the SLP.
We define the set $M$ of monomials as a union of $x_n$-chains:
\begin{align*}
  M = \{ v_i x_n^k \;;\; 1 \le i \le m, \; 0 \le k < p_i \}.
\end{align*}
This set $M$ is the standard monomials 
with respect to an ideal $I \in \cM_n(h)$ as follows.
It suffices to prove that 
any factor of monomials in $M$ again belongs to $M$.
For $u \in M$ not divisible by $x_n$,
any factor of $u$ is in $M$,
since $u$ is a standard monomial with respect to $\Ibar$.
For $u x_n^s \in M$ ($s \ge 1$),
the factor $u x_n^{s-1}$ sits in the same $x_n$-chain as $u x_n^s$,
and hence $u x_n^{s-1} \in M$ .
Let $v$ be a factor of $u$.
Then the factor $v x_n^s$ of $u x_n^s$ belongs to $M$,
since the $x_n$-chain starting with $v$ is not shorter than
that starting with $u$,
by the construction of $M$.
Thus we have proved that any factor of monomials in $M$ 
again belongs to $M$,
and this shows that the surjectivity of $\varphi$.

(ii)
In the proof above,
if $h$ is quasi-symmetric,
then there is only one choice of the lengths of $x_n$-chains
starting with $v_i$ and $v_j$,
where $\deg v_i = \deg v_j$.
Hence the map $\varphi$ is bijective.
\end{proof}
Next we restrict the mapping $\varphi$ in 
Lemma~\ref{lem:intersection-with-PPnm}
to the Borel-fixed ideals 
for which the quotient rings have the $k$-SLP.
We define the following collections of ideals in $R = \PPn$
for $0 \le k \le n$ and a given finite O-sequence $h$: 
\begin{align*}
  \begin{split}
    \cB_n^k(h) &=
    \{ I \in \cM_n(h) \;;\;
      \text{$I$ is Borel-fixed and $R/I$ has the $k$-SLP} \},
    \\
    \cA_n(h) &=
    \{ I \in \cM_n(h) \;;\;
      \text{$I$ is almost revlex} \},
  \end{split}
\end{align*}
where 0-SLP means the empty condition.
Clearly we have the following inclusions.
\begin{gather}
  \label{eq:inclusion-of-cSk-cBk-cA}
  \cB_n^0(h)  \supset \cdots \supset 
  \cB_n^{k-1}(h) \supset \cB_n^k(h) 
  \supset \cdots \supset \cB_n^n(h)  \supset  \cA_n(h) 
\end{gather}
Here the inclusion with $\cA_n$ is 
by Example~\ref{ex:almost-revlex-is-nSLP}.
Note that $\cB_n^k(h)$ is not empty in the case where
$h$ is a Hilbert function of some graded algebra with the $k$-SLP.

\begin{proposition}
\label{prop:kBorel-to-kBorel} 
Let $R = \PPn$ be the polynomial ring 
over a field of characteristic zero,
and $\Rbar = \PPnm$.
Let $k$ be an integer with $1 \le k \le n$.
Take an O-sequence $h = (h_0, h_1, \ldots, h_c)$,
which is a Hilbert function of some graded algebra with the $k$-SLP.
Define the map $\varphi_B$ as
\begin{align*}
  \varphi_B: \cB_n^k(h) \to \cB_{n-1}^{k-1}(\Delta h) 
  \qquad(\varphi_B: I \mapsto I \cap \Rbar). 
\end{align*}
Then we have the following:
\begin{itemize}
  \item[(i)]
  $\varphi_B$ is surjective,
  \item[(ii)]
  $\varphi_B$ is bijective, if $h$ is quasi-symmetric.
\end{itemize}
\end{proposition}
\begin{proof}
We first show that $\varphi_B$ is well-defined.
Take an ideal $I \in \cB_n^k(h)$.
It is clear that $I \cap \Rbar$ is Borel-fixed.
Since $x_n$ is a Lefschetz element of $R/I$
by Lemma~\ref{lem:wiebe-2.7}, 
it follows from Remark~\ref{rem:WLP-and-SLP}
that the Hilbert function of $\Rbar / (I \cap \Rbar)$ 
is equal to $\Delta h$.
Remark~\ref{rem:generic-elements-of-Borel-kSLP}
shows that $(R/I, x_n, \ldots, x_{n-k+2}, x_{n-k+1})$
has the $k$-SLP, if $R/I$ has the $k$-SLP.
Hence $\Rbar/(I \cap \Rbar)$ has the $(k-1)$-SLP,
and we have proved that $\varphi_B$ is well-defined.

Second, we prove the surjectivity.
Take an ideal $I \in \cB_{n-1}^{k-1}(h)$.
By Lemma~\ref{lem:intersection-with-PPnm} (i),
there is an inverse image $I$ of $\Ibar$,
whose quotient ring has the $k$-SLP.
We take the inverse image $I$ 
such that the standard monomials $v_i$ with respect to $\Ibar$ 
in the proof of  Lemma~\ref{lem:intersection-with-PPnm}
are ordered as $v_i <_{\text{revlex}} v_{i+1}$.
It remains to show that $I$ is Borel-fixed.
Let $v_i x_n^s$ ($s \ge 0$) be a standard monomial 
with respect to $I$, 
and suppose that $x_j | v_i x_n^s$.
To prove that $I$ is Borel-fixed,
we have to show that $(v_i/x_j) x_m  x_n^s$ is standard
for $j < m \le n$.
First, if $j < m < n$,  then $(v_i/x_j) x_m$ is a standard monomial
of $\Rbar$ with respect to $\Ibar$,
whose degree is equal to that of $v_i$ 
and satisfies that $(v_i/x_j) x_m <_{\text{revlex}} v_i$.
By the proof of Lemma~\ref{lem:intersection-with-PPnm},
the length of $x_n$-chain starting with $(v_i/x_j)x_m$
is not shorter than that starting with $v_i$.
Hence $(v_i/x_j)x_m x_n^s$ is standard.
Second, let $m = n$.
Note that $v_i/x_j$ is a standard monomial
of $\Rbar$ with respect to $\Ibar$,
whose degree is less than that of $v_i$. 
By the SL condition, 
the ending degree of $x_n$-chain starting with $v_i/x_j$
is not less than that starting with $v_i$.
Hence $(v_i/x_j) x_n^{s+1}$ is standard.
Thus we have proved that $I$ is Borel-fixed,
and $\varphi_B$ is surjective.

Third, the bijectivity of the case where $h$ is quasi-symmetric
follows from Lemma~\ref{lem:intersection-with-PPnm} (ii).
\end{proof}
\begin{remark}
\label{rem:3var-Borel-to-Borel}
Even when a given Hilbert function $h$ 
is not necessarily quasi-symmetric,
the map $\varphi_B$ of Proposition~\ref{prop:kBorel-to-kBorel}
happens to be bijective in the following cases.
When $n \le 2$, Borel-fixed ideals are unique 
for given Hilbert functions, 
and hence $\varphi_B$ is bijective.
Moreover, when $n=3$,
Borel-fixed ideals whose quotient rings have the SLP are unique 
by Theorem~\ref{thm:3var-uniqueness-of-Borel-SLP},
and hence $\varphi_B$ is bijective.
\end{remark}
\begin{remark}
\label{rem:almost-to-almost} 
Let $h$ be the Hilbert function of a quotient ring $R/I$
for some Artinian almost revlex ideal $I \subset R = \PPn$.
Then $\cA_n(h)$ contains only one ideal,
since almost revlex ideals are unique for given Hilbert functions.
The ideal $I \cap \PPnm$ is again almost revlex,
and hence $\cA_{n-1}(\Delta h)$ contains only one ideal.
Thus we have the following bijection $\varphi_A$.
\begin{align*}
  \begin{array}{cccl}
    \cA_n(h) & \overset{\varphi_A}{\to} & \cA_{n-1}(\Delta h) &
    \qquad(\varphi_A: I \mapsto I \cap \Rbar). 
  \end{array}
\end{align*}
\end{remark}

\subsection{Hilbert functions of graded algebras with the $k$-SLP}

We give a characterization of the Hilbert functions 
that can occur for graded $K$-algebras 
having the $k$-SLP or the $k$-WLP. 
Their characterizations are equal
as in the case of the SLP and the WLP
(Proposition~\ref{prop:hilbert-function-of-WLP-SLP}).
For a sequence $h = (h_0, h_1, \ldots, h_c)$ of positive integers,
define a sequence of the $t$-th difference $\Delta^t h$ by
\begin{align*}
  \Delta^t h &= \Delta(\Delta( \cdots \Delta(h) \cdots )) 
  \qquad \text{(apply $t$ times)},
\end{align*}
for a positive integer $t$.

\begin{proposition}
\label{prop:hilbert-fnct-of-nSLP}
Let $R = \PPn$ and $k$ be an integer with $1 \le k \le n$.
Let $h = (1, n, h_2, h_3, \ldots, h_c)$ be an O-sequence.
The following three conditions are equivalent:
\begin{itemize}
  \item[(i)]
  $h$ is a Hilbert function of some graded algebra with the $k$-SLP,
  \item[(ii)]
  $h$ is a Hilbert function of some graded algebra with the $k$-WLP,
  \item[(iii)]
  $h$ is a unimodal O-sequence,
  $\Delta^t h$ is a unimodal O-sequence for every integer $t$
  with $1 \le t < k$,
  and $\Delta^k h$ is an O-sequence.
\end{itemize}
\end{proposition}
\begin{proof}
It is clear that (i) implies (ii).
Using Remark~\ref{rem:WLP-and-SLP} and 
Proposition~\ref{prop:hilbert-function-of-WLP-SLP} repeatedly,
We can show that (ii) implies (iii).
We prove that (iii) implies (i).
Let $h$ be a sequence satisfying the condition of (iii).
Since we may take generic initial ideals,
there exists a Borel-fixed ideal $I'$ 
of the polynomial ring $R'$ in $(n-k)$ variables,
for which the Hilbert function of $R'/I'$ is equal to $\Delta^k h$.
Starting with $I'$, we can repeat taking one of inverse images
under the surjection $\varphi_B$ in Proposition~\ref{prop:kBorel-to-kBorel},
and we have a Borel-fixed ideal $I$ of $R$,
where $R/I$ has the $k$-SLP and has the Hilbert function $h$.
Thus we proved that (iii) implies (i).
\end{proof}
In addition, we have a characterization of the Hilbert functions
of quotient rings $R/I$ for Artinian almost revlex ideals $I$.
The characterization is the same as ideals with the $n$-WLP.
This result is an analogue of the result of 
Deery \cite{MR1381739} or
Marinari-Ramella \cite[Proposition 2.13]{MR1700543},
which gives the characterization of the Hilbert functions
for revlex ideals.
\begin{proposition}
\label{prop:hilbert-fnct-of-almost}
Let $R = \PPn$ and 
$h = (1, n, h_2, h_3, \ldots, h_c)$ an O-sequence.
The following four conditions are equivalent:
\begin{itemize}
  \item[(i)]
  $h$ is a Hilbert function of $R/I$ 
  for some almost revlex ideal $I$ of $R$,
  \item[(ii)]
  $h$ is a Hilbert function of some graded algebra with the $n$-SLP,
  \item[(iii)]
  $h$ is a Hilbert function of some graded algebra with the $n$-WLP,
  \item[(iv)]
  $h$ is a unimodal O-sequence,
  and $\Delta^k h$ is a unimodal O-sequence for every integer $k$
  with $1 \le k \le n$.
\end{itemize}
\end{proposition}
\begin{proof}
By Example~\ref{ex:almost-revlex-is-nSLP} and
Proposition~\ref{prop:hilbert-fnct-of-nSLP},
we have only to show that (ii) implies (i).
We prove it by induction on $n$.
The assertion is true for $n=1$.
Let $h$ be a Hilbert function of some graded algebra with the $n$-SLP.
By the assumption of induction,
we can take an almost revlex ideal $\Ibar$ of $\Rbar = \PPnm$,
where the Hilbert function of $\Rbar/\Ibar$ is equal to $\Delta h$.
Since $\Ibar$ is Borel-fixed 
and its quotient ring has the $(n-1)$-SLP,
there is an inverse image $I \subset R$
under the surjection $\varphi_B$
of Proposition~\ref{prop:kBorel-to-kBorel}.
In the construction of this inverse image, 
if we arrange $v_1, v_2, \ldots$ in the proof of
Lemma~\ref{lem:intersection-with-PPnm}
as $v_1 <_{\text{revlex}} v_2 <_{\text{revlex}} \cdots$,
then it follows from Proposition~\ref{prop:paraphrase-almost-revlex}
that $I$ is almost revlex.
Thus we have proved that (ii) implies (i).
\end{proof}

\begin{remark}
\label{rem:construct-n-SLP-hilbert}
The class of Hilbert functions of graded Artinian algebras
with the $n$-SLP
described in Proposition~\ref{prop:hilbert-fnct-of-almost}
is also characterized as follows using the condition (iv)
of Proposition~\ref{prop:hilbert-fnct-of-almost}:
\begin{itemize}
  \item[(i)]
  The finite sequences
  $(1,1,\ldots,1)$ form the class of Hilbert functions 
  of graded Artinian algebras $K[x_1]/I$ having the 1-SLP.
  \item[(ii)]
  All the Hilbert functions of $\PPn/I$ having the $n$-SLP
  can be obtained by
  applying the following manipulation $(n-1)$ times
  to sequences $(1,1,\ldots,1)$.
\end{itemize}
The manipulation is as follows:
\begin{quote}
  For a given sequence $a_0, a_1, a_2, \ldots, a_m$,
  make a finite sequence $b_j$ by
  \begin{gather*}
    b_0 = a_0, \quad
    b_1 = a_0+a_1, \quad\ldots, \quad
    b_m = a_0+\cdots+a_m,
    \\\text{and}~
    b_{m}, b_{m+1}, \cdots
    ~\text{are any weakly decreasing positive integers.}
  \end{gather*}
\end{quote}
\end{remark}

\subsection{Uniqueness of Borel-fixed ideals 
  and generic initial ideals in $\PPn$}

When $n \le 3$, we already know that Borel-fixed ideals of $\PPn$
whose quotient rings have the $n$-SLP
are the unique almost revlex ideals for given Hilbert functions.
Moreover, we have the following theorem for any $n$.
For a sequence $h$,
we use a convention that the 0-th difference $\Delta^0 h$ 
is $h$ itself.
\begin{theorem}
\label{thm:uniqueness-of-Borel-nSLP-symmetric} 
Let $I \subset R=\PPn$ be an Artinian Borel-fixed ideal 
whose quotient ring $R/I$ has the $n$-SLP,
and let $h$ be the Hilbert function of $R/I$. 
Suppose that the $k$-th difference $\Delta^k h$ is quasi-symmetric
for every integer $k$ with $0 \le k \le n-4$.
Then $I$ is the unique almost revlex ideal
for which the Hilbert function of $R/I$ is equal to $h$.
In particular, $I$ is determined only by the Hilbert function.
\end{theorem}
\begin{proof}
The theorem is already proved for $n \le 3$, and let $n \ge 4$.
Applying $\varphi_B$ of Proposition~\ref{prop:kBorel-to-kBorel}
to $I$ $(n-3)$-times,
we obtain the Borel-fixed ideal $I \cap \PPiii$,
whose quotient ring has the 3-SLP.
This Borel-fixed ideal is the unique almost revlex ideal 
for the Hilbert function $\Delta^{n-3} h$.
All $\varphi_B$ in this process are bijective
by Proposition~\ref{prop:kBorel-to-kBorel}.
For every $k$, note that $\Delta^k h$ is a Hilbert function 
for an algebra having the $(n-k)$-SLP,
and hence for the quotient ring by an almost revlex ideal
by Proposition~\ref{prop:hilbert-fnct-of-nSLP}.
Therefore all successive inverse images of $I \cap \PPiii$
are almost revlex by Remark~\ref{rem:almost-to-almost}.
Hence $I$ is almost revlex.
\end{proof}
In particular, 
we have the following uniqueness for Borel-fixed ideals
in the case of four variables.
\begin{corollary}
\label{cor:4var-uniqueness-of-symmetric-borel-2SLP}
Let $I \subset \PPiv$ be a Borel-fixed ideal,
for which $\PPiv/I$ has a quasi-symmetric Hilbert function $h$ ,
and has the 2-SLP.
Then $I$ is the unique almost revlex ideal 
for the Hilbert function $h$.
\qed
\end{corollary}
In Theorem~\ref{thm:uniqueness-of-Borel-nSLP-symmetric},
if we drop the condition for $\Delta^k h$ to be quasi-symmetric,
then the uniqueness does not necessarily hold as follows.
\begin{example}
\label{ex:4var-non-symmetric-Borel-2SLP-is-not-unique} 
(i)
There exist two different Borel-fixed ideals with the 4-SLP 
in $R=\PPiv$,
and their quotient rings have
the same non-quasi-symmetric Hilbert function.
Define the following ideals:
\begin{align*}
  I &= 
  (x_1^2, x_1x_2, x_2^3, x_2^2x_3, x_1x_3^2, x_2x_3^2, x_3^3,
    x_2^2x_4)
  +(x_1,x_2,x_3,x_4)^4, 
  \\
  J &=
  (x_1^2, x_1x_2, x_2^3, x_2^2x_3, x_1x_3^2, x_2x_3^2, x_3^3, 
    x_1x_3x_4)
  +(x_1,x_2,x_3,x_4)^4. 
\end{align*}
We can easily check that both $I$ and $J$ 
are Borel-fixed, have the 4-SLP,
and $R/I$ and $R/J$ have the same Hilbert function $h = (1,4,8,7)$.

Moreover this example shows that 
$\varphi_B$ of Proposition~\ref{prop:kBorel-to-kBorel}
is not injective for
$n = k = 4$ and  $h = (1,4,8,7)$,
since $I \cap \Rbar = J \cap \Rbar$ for $\Rbar = \PPiii$.

(ii)
The inverse images of $I$ and $J$ under $\varphi_B$ of 
Proposition~\ref{prop:kBorel-to-kBorel}
give two distinct Borel-fixed ideals,
whose quotient rings  have the 5-SLP,
and have the same symmetric Hilbert function
$h' = (1, 5, 13, 20, 13, 5, 1)$.
\end{example}
In the rest of this section,
we study generic initial ideals in the polynomial ring
$R = \PPn$ over a field of characteristic zero.
The following theorem,
which gives a uniqueness of generic initial ideals,
follows from
Theorem~\ref{thm:uniqueness-of-Borel-nSLP-symmetric} and 
Proposition~\ref{prop:I-is-kSLP-iff-Gin(I)-is-kSLP}.
\begin{theorem}
\label{thm:gin-of-symmetric-nSLP}
Let $I \subset R=\PPn$ be a graded Artinian ideal
whose quotient ring $R/I$ has the $n$-SLP,
and let $h$ be the Hilbert function of $R/I$. 
Suppose that the $k$-th difference $\Delta^k h$ is quasi-symmetric
for every integer $k$ with $0 \le k \le n-4$.
Then the generic initial ideal $\gin(I)$
with respect to the graded reverse lexicographic order 
is the unique almost revlex ideal for the Hilbert function $h$.
In particular, $\gin(I)$ is determined only by the Hilbert function.
\qed
\end{theorem}
In particular, we have the following corollary,
which corresponds to 
Corollary~\ref{cor:4var-uniqueness-of-symmetric-borel-2SLP}.
\begin{corollary}
\label{cor:4var-gin-of-symmetric-2SLP}
Let $I \subset \PPiv$ be a graded Artinian ideal 
whose quotient ring has the 2-SLP.
Suppose that the Hilbert function $h$ of $\PPiv/I$ is quasi-symmetric.
Then the generic initial ideal $\gin(I)$
with respect to the graded reverse lexicographic order 
is the unique almost revlex ideal for the Hilbert function $h$.
In particular, $\gin(I)$ is determined only by the Hilbert function.
\qed
\end{corollary}
\begin{remark}
We give a remark on Moreno's conjecture stated 
in \cite[Conjecture D]{math-0610649}.
The conjecture is as follows:
\begin{quote}
{\bf Moreno's Conjecture.}
If $f_1, f_2, \ldots, f_n \in R = \PPn$ is a generic sequence 
of homogeneous polynomials of given degree $d_1, d_2, \ldots, d_n$,
$I = (f_1, \ldots, f_n)$ and $J$ is the initial ideal of $I$
with respect to the graded reverse lexicographic order,
then $J$ is an almost revlex ideal.
\end{quote}
Cimpoea\c{s} \cite{math-0610649} proves that 
this conjecture is true for $n=3$.
When $n=4$,
if the quotient ring of any generic complete intersection
has the 2-SLP,
then the conjecture holds 
by Corollary~\ref{cor:4var-gin-of-symmetric-2SLP}.

When a socle degree is fixed,
Watanabe \cite{MR951211} shows that
the quotient rings of generic complete intersections in $\PPiii$
has the SLP.
In relation to Moreno's conjecture, 
the following question arises:
\begin{quote}
{\bf Question.}
Does any generic almost complete intersection in $\PPiii$
have the SLP?
\end{quote}
\end{remark}
\begin{example}
\label{ex:n-var-ci-SLP}
Let $R=\PPn$ be the polynomial ring 
over a field $K$ of characteristic zero. 
We consider complete intersections as follows. 
\begin{itemize}
\item[(a)] 
Let $f_1$ and $f_2$ be homogeneous polynomials of degree $d_i$ $(i=1,2)$, 
and let $g_3,\ldots,g_n$ be linear forms. 
Set $I=(f_1, f_2, f_3=g_3^{d_3}, \ldots, f_n=g_n^{d_n})$. 
Suppose that $\{f_1,f_2,g_3,\ldots,g_n\}$ is a regular sequence. 
Example 6.2 in \cite{harima-watanabe-jpaa} shows that $R/I$ has the SLP. 
\item[(b)] 
For $i=1,2,\ldots,n$, 
$f_i \in K[x_i,\ldots,x_n]$ be a homogeneous polynomial of degree $d_i$ 
which is a monic in $x_i$, 
and set $I=(f_1,f_2,\ldots,f_n)$. 
Then $R/I$ is always a complete intersection. 
Corollary 29 in \cite{MR2033004} 
and Corollary 2.1 in \cite{math-0506537} show that 
$R/I$ has the SLP. 
\end{itemize}
Now, let $k$ be an integer satisfying $1\leq k\leq n-2$ 
and suppose that 
$$
d_j\geq d_1+d_2+\cdots+d_{j-1}-(j-1)+1
$$
for all $j=n-k+1,n-k+2,\ldots,n$.  
Then we have the following. 
\begin{itemize}
\item[(i)]
$A=R/I$ has the $k$-SLP.  
\item[(ii)]
In particular, when $k=n-2$, 
$A$ has the $n$-SLP.
\item[(iii)]
The generic initial ideal of $I$ 
coincides with the unique almost revlex ideal 
determined by the Hilbert function of $A$. 
\item[(iv)] 
$\gin(x_1^{d_1}, x_2^{d_2}, \ldots, x_n^{d_n})=\gin(I)$. 
\end{itemize}
\end{example}
\begin{proof}
Choose a general linear form  $\overline{\ell} \in A$, 
where $\ell \in R$ is a linear from, 
satisfying the following conditions: 
\begin{itemize}
\item
$\{f_1, f_2, \ldots, f_{n-1}, \ell\}$ is a regular sequence, 
\item
$\overline{\ell}$ is a Lefschetz element of $A$. 
\end{itemize}
Set $S=R/\ell R$, and 
$$
B=A/\overline{\ell}A
=S/(\overline{f_1}, \overline{f_2}, 
\ldots, \overline{f_{n-1}}, \overline{f_n}). 
$$
Noting that 
$d_n\geq d_1+d_2+\cdots+d_{n-1}-(n-1)+1$ 
and $\{\overline{f_1}, \overline{f_2}, \ldots, \overline{f_{n-1}}\}$ 
is a regular sequence, 
it is easy to verify that 
$(\overline{f_1}, \overline{f_2}, \ldots, \overline{f_{n-1}},\overline{f_n})
=(\overline{f_1}, \overline{f_2}, \ldots, \overline{f_{n-1}})$
and 
$B=S/(\overline{f_1}, \overline{f_2}, \ldots, \overline{f_{n-1}})$ 
is a complete intersection. 
Hence $B$ also has the SLP by 
Example 6.2 in  \cite{harima-watanabe-jpaa},
Corollary 29 in \cite{MR2033004} 
and Corollary 2.1 in \cite{math-0506537}, 
and $A$ has the $2$-SLP. 
Repeating this argument, we have the assertion (i). 
The assertion (ii) follows from Remark~\ref{rem:k-SLP-k-WLP}. 
Next, we note that 
all $k$-th differences of the Hilbert function of $A$ are symmetric. 
Hence The assertion (iii) follows 
from Theorem~\ref{thm:gin-of-symmetric-nSLP}. 
The assertion (iv) is easy.   
\end{proof}
We conclude this section by an additional relation of initial ideals
with the $k$-WLP or the $k$-SLP.
Although this result is not used in the rest of this paper,
it is an analogue of 
Wiebe's result \cite[Proposition~2.9]{MR2111103}.
\begin{proposition}
\label{prop:I-is-kSLP-if-in(I)-is-kSLP}
Let $I$ be a graded Artinian ideal of $R=K[x_1,\ldots,x_n]$, 
let $\In(I)$ be the initial ideal of $I$ 
with respect to the graded reverse lexicographic order 
and let $1\leq k\leq n$. 
If $R/\In(I)$ has the $k$-WLP (resp. the $k$-SLP), 
then the same holds for $R/I$. 
\end{proposition}
\begin{proof}
We can give a proof using a similar idea to the proof of
\cite[Proposition~2.9]{MR2111103}.
First we consider the case of the $k$-WLP. 
Conca proves in
\cite[Theorem~1.1]{MR1948090}
that 
\begin{equation}\label{eq:conca1}
  \HH_{R/(I+(g_1,\ldots,g_i))}(t) \leq \HH_{R/(\In(I)+(g_1,\ldots,g_i))}(t)
\end{equation}
for generic linear forms $g_1,\ldots,g_i$ and all $t\geq 0$. 
In order to prove our claim, 
from Remark~\ref{rem:WLP-and-SLP} (i), 
it is enough to show that 
the Hilbert function of $R/(I+(g_1,\ldots,g_i))$ coincides with 
that of $R/(\In(I)+(g_1,\ldots,g_i))$ for every $i=1,2,\ldots,k$ 
under the assumption that the Hilbert function of 
$R/(I+(g_1,\ldots,g_{i-1}))$ is equal to that of 
$R/(\In(I)+(g_1,\ldots,g_{i-1}))$. 
Set $h = \HH_{R/(I+(g_1,\ldots,g_{i-1}))}$.
Using an argument of induction, 
we may assume that $R/(\In(I)+(g_1,\ldots,g_{i-1}))$ has the WLP.
Hence it follows from Remark~\ref{rem:WLP-and-SLP} (i) that 
the Hilbert function of $R/(\In(I)+(g_1,\ldots,g_i))$ 
is equal to $\Delta h$.
Furthermore one can easily check that 
\begin{align*}
  (\Delta h)_t \leq \HH_{R/(I+(g_1,\ldots,g_i))}(t)
\end{align*}
for all $t \ge 0$.
Hence it follows from (\ref{eq:conca1}) that 
\begin{align*}
  \HH_{R/(I+(g_1,\ldots,g_i))}(t) = \HH_{R/(\In(I)+(g_1,\ldots,g_i))}(t)
\end{align*}
for all $t\geq 0$. 

Next we consider the case of the $k$-SLP. 
From the proof of \cite[Proposition~2.9]{MR2111103}, 
we have  
\begin{align}
  \label{eq:conca2}
  \HH_{R/(I+(g_1,\ldots,g_{i-1},g_i^s))}(t) \leq 
  \HH_{R/(\In(I)+(g_1,\ldots,g_{i-1},g_i^s))}(t)
\end{align}
for generic linear forms $g_1,\ldots,g_i$, 
$s \geq 1$ and all $t\geq 0$. 
In order to prove our claim, 
it is enough to show from Remark~\ref{rem:WLP-and-SLP} (ii) that 
the Hilbert function of $R/(I+(g_1,\ldots,g_{i-1},g_i^s))$ 
coincides with that of $R/(\In(I)+(g_1,\ldots,g_{i-1},g_i^s))$ 
for every $i=1,2,\ldots,k$ 
under the assumption that the Hilbert function of
 $R/(I+(g_1,\ldots,g_{i-1}))$ 
is equal to that of $R/(\In(I)+(g_1,\ldots,g_{i-1}))$. 
Set $h = \HH_{R/(I+(g_1,\ldots,g_{i-1}))}$.
Using an argument of induction, 
we may assume that $R/(\In(I)+(g_1,\ldots,g_{i-1}))$ has the SLP.  
Hence, it follows from Remark~\ref{rem:WLP-and-SLP} (ii) that 
the Hilbert function of $R/(\In(I)+(g_1,\ldots,g_{i-1},g_i^s))$ 
is equal to the sequence $(b_t)_{t \ge 0}$:
\begin{align*}
  b_t = \max\{ h_t - h_{t-s}, 0 \},
\end{align*}
where $h_t = 0$ for $t < 0$.
Furthermore one can easily check that 
\begin{align*}
  b_t \leq \HH_{R/(I+(g_1,\ldots,g_{i-1},g_i^s))}(t)
\end{align*}
for all $t \ge 0$.
Hence it follows from (\ref{eq:conca2}) that 
\begin{align*}
  \HH_{R/(I+(g_1,\ldots,g_{i-1},g_i^s))}(t) 
  = \HH_{R/(\In(I)+(g_1,\ldots,g_{i-1},g_i^s))}(t)
\end{align*}
for all $t \geq 0$. 
\end{proof}

\section{An extremal property of graded Betti numbers and the $k$-WLP}

In the rest of this paper,
we study graded Betti numbers for monomial ideals.
The goal is Theorem~\ref{thm-bounds}
on the maximality of graded Betti numbers.
We give a sharp upper bound on the graded Betti numbers 
of graded Artinian algebras 
with the $k$-WLP and a fixed Hilbert function.
The upper bounds are achieved by the quotient rings 
by Borel-fixed ideals having the $k$-SLP.
In particular, when $k=n$,
almost revlex ideals give the upper bounds.
\subsection{Graded Betti numbers of stable ideals and the $k$-WLP}

Let $R = \PPn$ be the polynomial ring over a field 
of characteristic zero,
and $\frakm = (x_1, x_2, \ldots, x_n)$ be the maximal graded ideal.
We recall $\frakm$-full ideals,
which are first defined by Watanabe in \cite{MR894414},
and studied in \cite{MR1075117}.
$\frakm$-Full ideals are useful in studying graded Betti numbers.
\begin{definition}[\cite{MR894414}]
\label{defn:m-full}
An ideal $I \subset R$ is called an $\frakm$-{\it full ideal},
if there exists an element $z \in R$ such that
$\frakm I : z = I$.

It is known that $z$ can be taken as a general homogeneous element 
of degree one, if $I$ is graded.
Note that an inclusion $\frakm I : z \supset I$ always holds,
if $z$ does not have a constant term.
\end{definition}
\begin{example}
\label{ex:m-full}
(i)
An easy example of $\frakm$-full ideals is a prime ideal $\frakp$
\cite[Remark 1]{MR894414}.
Indeed, $\frakp$ is an associated prime of $\frakm\frakp$,
and hence there is $z$ such that $\frakm\frakp:z = \frakp$.

(ii)
If $I$ is a stable ideal,
then $I$ is $\frakm$-full as follows:
Set $z=x_n$ and let $m \in \frakm I:z$ be a monomial.
Since $z m \in \frakm I$, we have $z m=x_i m^\prime$ 
for some $x_i$ and a monomial $m' \in I$.
Hence $m=x_i m' / z \in I$, since $I$ is stable.
Thus $I$ is $\frakm$-full.

(iii)
If $I \subset R$ is a monomial ideal and $(R/I, x_n)$ has the SLP,
then $I$ is $\frakm$-full as follows:
Set $z = x_n$,
and let $m$ be a monomial in $\frakm I : x_n$.
We prove that $m \in I$.
There exist a minimal generator $m'$ of $I$ 
and a monomial $u \in R$ with $\deg u \ge 1$,
such that $x_n m = um'$.
If $x_n$ divides $u$, then $m = (u/x_n)m'$ belongs to $I$.
If $x_n$ does not divide $u$,
then $m'/x_n = m/u$ has the same exponent of $x_n$ as that of $m$,
and satisfies that $\deg m' \le \deg m$.
Assume that $m \not\in I$,
and consider the following two $x_n$-chains:
(a) the $x_n$-chain containing $m$,
and (b) the $x_n$-chain ending in $m'/x_n$.
Then the ending degree of (a) is greater than that of (b)
by $\deg m' \le \deg m$,
and the starting degree of (a) is greater than that of (b),
since the exponent of $x_n$ in $m$ and that in $m'/x_n$ are the same,
and $\deg (m'/x_n) < \deg m$.
This contradicts the SL condition, and hence $m \in I$.
\end{example}
J. Watanabe proved formulas for Betti numbers
of $\frakm$-full ideals \cite[Corollary~8, 9]{MR1075117}.
It should be remarked that these formulas are applicable
even when the ideal is not graded.
For graded $\frakm$-full ideals,
we obtain a formula for graded Betti numbers as follows.
We also remark that this formula recovers
the formula of Eliahou and Kervaire \cite{MR1037391}
for graded Betti numbers for stable ideals.
\begin{proposition}[{\cite[Corollary~8]{MR1075117}}]
\label{prop:betti-of-m-full-is-reduced-to-low-rank}
Let $I$ be a graded $\frakm$-full ideal of $R$,
and $z \in R_1$ an element satisfying $\frakm I: z = I$.
Set $\Rbar = R/z R$, and $\Ibar = (I+z R)/z R$.
Write $I = (f_1, f_2, \ldots, f_r, z f_{r+1}, \ldots, z f_s)$
by homogeneous polynomials,
where they are minimal generators of $I$,
and first $r$ polynomials are minimal generators of $\Ibar$.
Then we have the following formula:
\begin{align}
  \label{eq:betti-of-m-full}
\begin{split}
  \beta_{i,i+j}(R/I) &=
  \beta_{i,i+j}(\Rbar/\Ibar) +
  \binom{n-1}{i-1} \times c_{j+1},
  \\
  c_j &= \#\{ t \;;\; r <  t \le s, \deg(z f_t) = j \},
\end{split}
\end{align}
where
$\beta_{i,i+j}(R/I) = \dim_K [\Tor_i^R(R/I,K)]_{i+j}$,
and
$\beta_{i,i+j}(\Rbar/\Ibar) = 
\dim_K [\Tor_i^{\Rbar}(\Rbar/\Ibar,K)]_{i+j}$.

In particular, if $I$ is a monomial ideal and $z = x_n$,
then we have
\begin{align}
  \label{eq:cj-in-betti-of-m-full}
  c_j = \#\{ \text{the minimal generators of degree $j$
      divisible by $x_n$} \}.
\end{align}
\qed
\end{proposition}
The following lemma shows a relation of $\frakm$-full ideals
with $x_n$-chains.
In particular, stable ideals have the property of this lemma.
Although most results in this section concerning stable ideals
can be proved for monomial $\frakm$-full ideals satisfying
$\frakm I : x_n = I$,
we state such results only for stable ideals 
for simplicity.
\begin{lemma}
\label{lem:m-full-is-xn-chain} 
Let $I \subset R = \PPn$ be a monomial Artinian $\frakm$-full ideal
satisfying $\frakm I : x_n = I$,
and $u, u x_n, \ldots, u x_n^{s-1}$ an $x_n$-chain
with respect to $I$,
where $u$ is a monomial not divisible by $x_n$, and $s \ge 1$.
Then $u x_n^s$ is a member of the minimal generators of $I$.
\end{lemma}
\begin{proof}
Let $I$ be a monomial $\frakm$-full ideal, and $\frakm I : x_n = I$.
Let $u \in R$ be a monomial not divisible by $x_n$,
and suppose that $u x_n ^{s-1} \not\in I$ and $u x_n^s \in I$.
There exists a minimal generator $v x_n^t$ dividing $u x_n^s$,
and $v$ is not divisible by $x_n$.
Namely, $v | u$ and $t \le s$.
We have $t = s$, since $t < s$ contradicts that
$u x_n^{s-1} \not\in I$.
Assume that $\deg v < \deg u$.
Then there exists $i<n$ such that $x_i v | u$,
and hence $x_i v x_n^{s-1} | u x_n^{s-1}$.
Since $x_i v x_n^s \in \frakm I$ and $I$ is $\frakm$-full, 
we have $x_i v x_n^{s-1} \in I$,
and hence $u x_n^{s-1} \in I$.
This is a contradiction, and we have $v = u$.
This means that $u x_n^s$ is a minimal generator.
\end{proof}
\begin{remark}
\label{rem:xn-chain-ideal-with-WLP}
Let $I$ be an Artinian stable ideal of $R = \PPn$,
and suppose that $(R/I, x_n)$ has the WLP.
Let $h = (h_0, h_1, \ldots, h_c)$ be the Hilbert function of $R/I$.
It follows from Example~\ref{ex:m-full} (ii),
Lemma~\ref{lem:m-full-is-xn-chain} and the WL condition that
\begin{itemize}
  \item[(i)]
  the number of the minimal generators of degree $d$
  divisible by $x_n$
  is equal to $\max\{ h_{d-1}-h_d, 0\}$.
\end{itemize}
We can also say that
\begin{itemize}
  \item[(ii)]
  the minimal generators divisible by $x_n$ occur only at degrees
  $u_2, u_3, \ldots$ defined in Equation~(\ref{eq:def-of-ui}),
  and the number of such generators of degree $u_j$ ($j \ge 2$)
  is equal to $h_{u_j-1}-h_{u_j}$.
\end{itemize}
In particular, the degrees of minimal generators divisible by $x_n$
are determined only by the Hilbert function.
\end{remark}
We have the following proposition
by Proposition~\ref{prop:betti-of-m-full-is-reduced-to-low-rank}.
It gives graded Betti numbers for stable ideals
whose quotient rings have the WLP.
\begin{proposition}
\label{prop:betti-of-stable-WLP-is-reduced-to-low-rank}
Let $I \subset R$ be an Artinian stable ideal,
for which $R/I$ has the Hilbert function $h = (h_0, h_1, \ldots, h_c)$,
and $(R/I, x_n)$ have the WLP.
Let $\Rbar = \PPnm$ and $\Ibar = I \cap \Rbar$.
We have the following.
\begin{itemize}
  \item[(i)]
  The graded Betti numbers $\beta_{i,i+j}(R/I)$ of $R/I$
  is given as follows:
  \begin{align*}
    \beta_{i,i+j}(R/I) &= 
    \beta_{i,i+j}(\Rbar/\Ibar) + \binom{n-1}{i-1} \times c_{j+1}
    \qquad (i,j \ge 0),
    \\
    c_j &= \max\{h_{j-1}-h_j, 0\},
  \end{align*}
  where we use the convention that $h_{-1} = 0$.
  \item[(ii)]
  By the same $c_j$, the last graded Betti numbers $\beta_{n,n+j}(R/I)$
  is given as follows:
  \begin{align*}
    \beta_{n,n+j}(R/I) = c_{j+1} \qquad (j \ge 0).
  \end{align*}
  In particular, they are determined only by the Hilbert function.
\end{itemize}
\end{proposition}
\begin{proof}
Let $I \subset R$ be an Artinian stable ideal,
and suppose that $(R/I, x_n)$ has the WLP.
In this case, the degrees of minimal generators divisible by $x_n$
are determined only by the Hilbert function $h$ of $R/I$
thanks to Remark~\ref{rem:xn-chain-ideal-with-WLP}.
Namely the constant $c_j$ 
in Proposition~\ref{prop:betti-of-m-full-is-reduced-to-low-rank}
is written as $c_j = \max\{h_{j-1}-h_j, 0\}$.
We thus obtain (i).
We have (ii), since $\beta_{n,n+j}(\Rbar/\Ibar) = 0$.
\end{proof}
We easily generalize
Proposition~\ref{prop:betti-of-stable-WLP-is-reduced-to-low-rank}
to the case of the $k$-WLP.
\begin{notation}
For a unimodal O-sequence $h=(h_0,h_1,\ldots,h_c)$, 
we define 
$$
c^{(h)}_j=\max\{h_{j-1}-h_j,0\} 
$$
for all $j=0,1,\ldots,c$, 
where $h_{-1}=0$. 
\end{notation}

\begin{proposition}\label{Betti-of-kWLP}
Let $I$ be an Artinian 
Borel-fixed ideal of $R=K[x_1,x_2,\ldots,x_n]$ 
and suppose that $R/I$ has the $k$-WLP. 
\begin{itemize}
\item[(i)]
Let $k<n$. 
Set $R^\prime=K[x_1,x_2,\ldots,x_{n-k}]$ and $I^\prime=I\cap R^\prime$. 
Then 
\begin{align*}
\beta_{i,i+j}(R/I)  =  
\beta_{i,i+j}(R^\prime/I^\prime)+
\tbinom{n-k}{i-1} 
\cdot c^{(\Delta^{k-1}h)}_{j+1} + 
\cdots+
\tbinom{n-2}{i-1}
\cdot c^{(\Delta h)}_{j+1}
+
\tbinom{n-1}{i-1}
\cdot c^{(h)}_{j+1}.
\end{align*}
\item[(ii)]
Let $k=n$. Then we have 
\begin{align*}
\beta_{i,i+j}(R/I)  =  
\tbinom{0}{i-1}
\cdot c^{(\Delta^{n-1}h)}_{j+1}
+
\tbinom{1}{i-1}
\cdot c^{(\Delta^{n-2}h)}_{j+1}
+\cdots+
\tbinom{n-2}{i-1}
\cdot c^{(\Delta h)}_{j+1}
+
\tbinom{n-1}{i-1}
\cdot c^{(h)}_{j+1}.
\end{align*}
In particular, $\beta_{i,i+j}(R/I)$ is determined 
only by the Hilbert function.
\end{itemize}

\end{proposition}
\begin{proof}
Note that $(R/I, x_n, x_{n-1}, \ldots, x_{n-k+1})$ has the $k$-WLP
by Remark~\ref{rem:generic-elements-of-Borel-kSLP}.
Repeating the same argument as the proof of 
Proposition~\ref{prop:betti-of-stable-WLP-is-reduced-to-low-rank}, 
we can get these equalities. 
\end{proof}
\begin{corollary}
Let $I$ be a graded Artinian ideal of $R=K[x_1,x_2,\ldots,x_n]$ 
and suppose that $R/I$ has the $k$-WLP. 
\begin{itemize}
\item[(i)]
Let $k<n$. 
Set $R^\prime=K[x_1,x_2,\ldots,x_{n-k}]$ 
and $I^\prime=\gin(I)\cap R^\prime$. 
Then 
\begin{align*}
\beta_{i,i+j}(R/\gin(I)) & = & 
\beta_{i,i+j}(R^\prime/I^\prime)+
\tbinom{n-k}{i-1} 
\cdot c^{(\Delta^{k-1}h)}_{j+1} + 
\cdots+
\tbinom{n-2}{i-1}
\cdot c^{(\Delta h)}_{j+1}
+
\tbinom{n-1}{i-1}
\cdot c^{(h)}_{j+1}.
\end{align*}
\item[(ii)]
Let $k=n$. Then we have 
\begin{align*}
\beta_{i,i+j}(R/\gin(I))  =  
\tbinom{0}{i-1}
\cdot c^{(\Delta^{n-1}h)}_{j+1}
+
\tbinom{1}{i-1}
\cdot c^{(\Delta^{n-2}h)}_{j+1}
+\cdots+
\tbinom{n-2}{i-1}
\cdot c^{(\Delta h)}_{j+1}
+
\tbinom{n-1}{i-1}
\cdot c^{(h)}_{j+1}.
\end{align*}
In particular, $\beta_{i,i+j}(R/\gin(I))$ is determined 
only by the Hilbert function.
\end{itemize}
\end{corollary}

\begin{proof}
This follows from Propositions~\ref{prop:I-is-kSLP-iff-Gin(I)-is-kSLP}
and \ref{Betti-of-kWLP}. 
\end{proof}

\subsection{Maximality of graded Betti numbers and the $k$-WLP}

\begin{notation and remark}
Let $h$ be the Hilbert function of a graded Artinian
$K$-algebra $R/I$.  Then there is the uniquely
determined lex-segment ideal $J \subset R$ such that
$R/J$ has
${h}$ as its Hilbert function.  We define
$$
\beta_{i,i+j} ({h},R) = \beta_{i,i+j} (R/J).
$$
The numbers $\beta_{i,i+j} ({h},R)$ can be computed numerically
without considering lex-segment ideals.  Explicit formulas can be found in
\cite{MR1037391}.
\end{notation and remark}
We give a sharp upper bound on the Betti numbers 
among graded Artinian $K$-algebras having the $k$-WLP. 
Moreover the upper bound is achieved by a graded Artinian $K$-algebra
with the $k$-SLP.
For $k=1$, this theorem was first proved by 
\cite[Theorem 3.20]{MR1970804}.
\begin{theorem}\label{thm-bounds} 
\begin{itemize}
  \item[(i)] 
  Let $A=R/I$ be a graded Artinian $K$-algebra with the $k$-WLP 
  and put $R^\prime=K[x_1,\ldots,x_{n-k}]$. 
  Then the graded Betti numbers of $A$ satisfy 
  \begin{align} 
    \label{eq:thm-bounds-1}
    \begin{split}
      \beta_{i,i+j}(A) 
      & \le
	\beta_{i,i+j}(\Delta^{k}h, R^\prime)+
	\tbinom{n-k}{i-1} 
	\cdot c^{(\Delta^{k-1}h)}_{j+1} 
	+\cdots 
	\\[1ex] & \qquad
	+
	\tbinom{n-2}{i-1}
	\cdot c^{(\Delta h)}_{j+1}
	+
	\tbinom{n-1}{i-1}
	\cdot c^{(h)}_{j+1}   \mbox{\hspace{8ex}if $k<n$,}
      \end{split}
      \\ 
      \intertext{and}
      \label{eq:thm-bounds-2}
      \begin{split}
      \beta_{i,i+j}(A)  
      & \le
	\tbinom{0}{i-1}
	\cdot c^{(\Delta^{n-1}h)}_{j+1}
	+
	\tbinom{1}{i-1}
	\cdot c^{(\Delta^{n-2}h)}_{j+1}
	+\cdots
	\\[1ex] & \qquad
	+
	\tbinom{n-2}{i-1}
	\cdot c^{(\Delta h)}_{j+1}
	+
	\tbinom{n-1}{i-1}
	\cdot c^{(h)}_{j+1}   \mbox{\hspace{8ex} if $k=n$}. 
      \end{split}
  \end{align}

  \item[(ii)]
  Let $h$ be an O-sequence such that 
  there is a graded Artinian $K$-algebra $R/J$ having the $k$-WLP  
  and $h$ as Hilbert function. 
  Then there is a Borel-fixed ideal $I$ of $R$ such that 
  $R/I$ has the $k$-SLP, 
  the Hilbert function of $R/I$ is $h$ 
  and the equality holds in (i) for all integers $i,j$. 
\end{itemize}
\end{theorem}

\begin{proof}
(i) 
Since $R/I$ has the $k$-WLP, 
it follows by Proposition~\ref{prop:I-is-kSLP-iff-Gin(I)-is-kSLP} that 
\\ 
$(R/\gin(I), x_n,\ldots,x_{n-k+1})$ 
has the $k$-WLP. 
Put $I^\prime=\gin(I)\cap R^\prime$. 
Noting that 
the Hilbert function of $R^\prime/I^\prime$ 
coincides with $\Delta^k h$, 
we have by \cite{MR1218500}, \cite{MR1218501}, \cite{MR1415019} that
$$
\beta_{i, i+j}(R^\prime/I^\prime) \leq 
\beta_{i,i+j}(\Delta^{k}h, R^\prime)
$$
for all $i,j$. 
Hence Proposition \ref{Betti-of-kWLP} shows that 
the Betti numbers of $R/\gin(I)$ satisfy the inequalities in (i). 
Thus our claim follows from the well known fact 
\cite{MR1648665},
$\beta_{i,i+j}(R/I)\leq\beta_{i,i+j}(R/\gin(I))$ 
for all $i,j$.  

(ii) 
We first consider the case where $k=n$. 
Let $I$ be the unique almost revlex ideal of $R$ 
whose quotient ring has the same Hilbert function $h$. 
Then $I$ is the desired ideal 
from Example~\ref{ex:almost-revlex-is-nSLP}
and Proposition \ref{Betti-of-kWLP} (ii). 

Let $k<n$.
Let  $I^{(k)}$ be the unique lex-segment ideal of $R'$,
where the Hilbert function of $R'/I^{(k)}$ is equal to $\Delta^k h$.
We can take an inverse image $I^{(k-1)}$ of $I^{(k)}$ 
under $\varphi_B$ of Proposition~\ref{prop:kBorel-to-kBorel}.
Then $I^{(k-1)}$ is an ideal of $K[x_1,x_2,\ldots, x_{n-k+1}]$,
and its quotient ring has the SLP
and $\Delta^{k-1} h$ as the Hilbert function.
Repeating the same procedure $k$ times,
we obtain a Borel-fixed ideal $I^{(0)}$ of $R$,
whose quotient ring has the $k$-SLP
and $h$ as the Hilbert function.
This ideal gives the upper bound of the graded Betti numbers
as desired.
\end{proof}

The following is an immediate consequence of Theorem \ref{thm-bounds}.  
\medskip

\begin{corollary}
\label{cor:maximal-Betti-number}
Let $h$ be the Hilbert function of  
a graded Artinian $K$-algebra $R/J$ having the $n$-WLP (resp. $n$-SLP). 
Let $I$ be the unique almost revlex ideal of $R$ 
whose quotient ring has the same Hilbert function $h$. 
Then $R/I$ has the maximal Betti numbers 
among graded Artinian $K$-algebras with the same Hilbert function $h$ 
and the $n$-WLP (resp. $n$-SLP).
\qed
\end{corollary}
\begin{corollary}
Let $I$ be a graded Artinian ideal of $R=K[x_1,x_2,\ldots,x_n]$ 
and suppose that $R/I$ has the $n$-WLP (resp. $n$-SLP). 
Then $R/\gin(I)$ has the maximal Betti numbers 
among graded Artinian $K$-algebras with the same Hilbert function $R/I$ 
and the $n$-WLP (resp. $n$-SLP).
\end{corollary}
\begin{proof}
Proposition~\ref{prop:I-is-kSLP-iff-Gin(I)-is-kSLP},
Proposition~\ref{Betti-of-kWLP} 
and Theorem~\ref{thm-bounds}
prove this corollary.
\end{proof}

\bibliographystyle{alpha}
\bibliography{math}

\end{document}